\documentclass[a4paper,12pt]{article}
\usepackage{amsmath}
\usepackage{amsfonts}
\usepackage{amssymb}
\usepackage{latexsym}
\usepackage{epsfig}
\usepackage{graphicx}
\usepackage{oldgerm}
\usepackage{theorem}
\usepackage{bm}

\def\N{\mathbb{N}}
\def\C{\mathbb{C}}
\def\R{\mathbb{R}}

\def\E{\mathbb{E}}
\def\P{\mathbb{P}}

\def\sm{\mathsf{m}}
\def\su{\mathsf{u}}
\def\sh{\mathsf{h}}
\def\sw{\mathsf{w}}
\def\so{\mathsf{o}}
\def\sa{\mathsf{a}}
\def\sx{\mathsf{x}}
\def\s{\mathsf{s}}

\def\bra{\langle}
\def\ket{\rangle}

\def\bLambda{\bm{\Lambda}}
\def\blambda{\bm{\lambda}}
\def\bkappa{\boldsymbol{\kappa}}
\def\1{{\bf 1}}

\def\L{\bm{L}}
\def\bR{\bm{R}}
\def\u{\bm{u}}
\def\v{\bm{v}}
\def\br{\bm{r}}
\def\bl{\bm{l}}

\def\cO{\mathcal{O}}
\def\cM{\mathcal{M}}
\def\cN{\mathcal{N}}
\def\cB{\mathcal{B}}

\def\cC{\mathcal{C}}

\def\cF{\mathcal{F}}
\def\cM{\mathcal{M}}

\def\rR{{\rm R}}
\def\rI{{\rm I}}
\def\Tr{{\rm Tr} \,}
\def\diag{{\rm diag}}
\def\rO{{\rm O}}

\def\supp{{\rm supp} \,}

\def\law={\stackrel{\rm (law)}{=}}

\def\blambdabar{\overline{\blambda}}
\def\zbar{\overline{z}}
\def\wbar{\overline{w}}

\def\zetabar{\overline{\zeta}}

\theorembodyfont{\itshape}

\newtheorem{thm}{Theorem}[section]
\newtheorem{lem}[thm]{Lemma}
\newtheorem{cor}[thm]{Corollary}
\newtheorem{prop}[thm]{Proposition}

\newtheorem{df}[thm]{Definition}

\newcommand{\SSC}[1]{\section{#1}\setcounter{equation}{0}}
\newcommand{\qed}{\hbox{\rule[-2pt]{3pt}{6pt}}}

\newfont{\bg}{cmr10 scaled\magstep4}

\newcommand{\bigzerou}{\smash{\lower1.7ex\hbox{\bg 0}}}

\begin{document}

\title{\Large{\bf
Eigenvalues, eigenvector-overlaps, and \\
regularized Fuglede--Kadison determinant of \\
the non-Hermitian matrix-valued Brownian motion
}
}
\author{
Syota Esaki
\footnote{
Mathematical Sciences Program, 
Department of Science and Technology, 
Faculty of Science and Technology,  
Oita University, Oita 870-1192, Japan;
e-mail: sesaki@oita-u.ac.jp 
}
\quad
Makoto Katori
\footnote{
Department of Physics,
Faculty of Fundamental Science and Engineering,
School of Science and Engineering,
Chuo University, 
Kasuga, Bunkyo-ku, Tokyo 112-8551, Japan; 
e-mail: makoto.katori.mathphys@gmail.com}
\quad
Satoshi Yabuoku
\footnote{
Department of Applied Mathematics, 
Fukuoka University, Fukuoka 814-0180, Japan; 
e-mail:
yabuoku@fukuoka-u.ac.jp
}
}

\date{5 April 2026}
\pagestyle{plain}
\maketitle

\begin{abstract}
The non-Hermitian matrix-valued Brownian motion
is the stochastic process of a random matrix
whose entries are given by independent complex Brownian motions.
The bi-orthogonality relation is imposed between
the right and the left eigenvector processes, which allows for
their scale transformations with an invariant eigenvalue process.
The eigenvector-overlap process is a 
Hermitian matrix-valued process,
each element of which is given by 
a product of an overlap of right eigenvectors
and that of left eigenvectors.
We derive a set of stochastic differential equations (SDEs)
for the coupled system of the eigenvalue process and the
eigenvector-overlap process and prove 
the scale-transformation invariance of the obtained SDE system.
The Fuglede--Kadison (FK) determinant associated with the present
matrix-valued process is regularized by introducing 
an auxiliary complex variable. This variable is necessary
to give the stochastic partial differential equations (SPDEs) 
for the time-dependent random field defined by
the regularized FK determinant 
and for its squared and logarithmic variations.
Time-dependent point process of eigenvalues
and its variation weighted by the diagonal elements of
the eigenvector-overlap process are 
related to the derivatives of the
logarithmic regularized FK-determinant random-field.
We also discuss the PDEs obtained by averaging the SPDEs. 

\vskip 0.2cm



\end{abstract}
\vspace{3mm}
\normalsize

\SSC
{Introduction and Results}
\label{sec:introduction}
One of the recent topics of
random matrix theory \cite{For10,Meh04} 
is \textit{random matrix dynamics}, 
which has been extensively
studied in statistical mechanics 
and probability theory.
In the \textit{Gaussian unitary ensemble} for
Hermitian random matrices,
algebraically independent elements of a matrix
are following independent Gaussian distributions.  
The system of eigenvalues of such a
Hermitian matrix-valued Gaussian random variable 
is identified with the \textit{one-dimensional log-gas},
since the eigenvalues of Hermitian matrix are
distributed on the real line and their repulsive
interactions are given by the logarithmic potential.
If we replace the independent Gaussian random variables
by the independent Brownian motions (BMs) at
any elements of a Hermitian matrix,
then we have a Hermitian matrix-valued BM.
Its eigenvalue process is described by a system of
stochastic differential equations (SDEs)
called the \textit{Dyson BM} \cite{Dys62,Kat_Springer}.
It is considered as the \textit{one-dimensional dynamical 
log-gas}, since the drift terms are given by
derivatives of log-potential.
The \textit{complex Ginibre ensemble} 
of non-Hermitian matrix-valued Gaussian 
random variables \cite{BF24,Gin65}
provides eigenvalue distributions
on the complex plane, which are regarded as
\textit{two-dimensional log-gases}.
They are also called the 
\textit{planar Coulomb gases},
since the two-dimensional Coulomb potential
of electrostatics is logarithmic.
Hence, it is natural and fruitful
to study the systems of SDEs
representing the two-dimensional Coulomb gases
\cite{BCF18,Osa12,Osa13}.
We claim that, however, such SDE systems
are not random matrix dynamics 
for non-Hermitian systems in the
sense of the Dyson BM for Hermitian systems. 
In the present paper, we start from the 
non-Hermitian matrix-valued BM
and study its eigenvalue process
and the associated processes. 
We will develop the previous studies
of random matrix dynamics reported by 
\cite{BGNTW16,BD20,Burda14,Burda15,GW18,Yab20} 
in this direction, and clarify that 
the non-Hermitian dynamics is quite different
from the Hermitian dynamics.

\subsection{
The non-Hermitian matrix-valued BM
and associated processes}
\label{sec:processes}

Let $N \in \N:=\{1,2, \dots\}$. 
Consider $2N^2$ independent one-dimensional standard
BMs,
$(B^{\rR}_{jk}(t))_{t \geq 0}$,
$(B^{\rI}_{jk}(t))_{t \geq 0}$, $1 \leq j, k \leq N$ defined in 
a filtered probability space $(\Omega, \cF, \{\cF_t\}_{t \geq 0}, \P)$,
where the expectation with respect to $\P$ is denoted by $\E$.
Let $i:=\sqrt{-1}$ and we define the 
\textit{$N \times N$ non-Hermitian
matrix-valued BM} by
\begin{equation}
M(t)=(M_{jk}(t))_{1 \leq j, k \leq N}
:= \left( \frac{1}{\sqrt{2N}}
(B^{\rR}_{jk}(t)+ i B^{\rI}_{jk}(t)) \right)_{1 \leq j, k \leq N},
\quad t \geq 0, 
\label{eq:Mt}
\end{equation}
which starts from a deterministic matrix $M(0)=(M_{jk}(0))_{1 \leq j, k \leq N}$.
We write the increment of $M(t)$ as
$
dM(t)=(dM_{jk}(t))_{1 \leq j, k \leq N}
=((dB^{\rR}_{jk}(t)+ i dB^{\rI}_{jk})/\sqrt{2N})_{1 \leq j, k \leq N}. 
$
By definition, the cross-variations of increments of 
the elements of $M(t)$ are given by
\begin{align}
\bra d M_{jk}, dM_{\ell m} \ket_t
&=0, 
\label{eq:variation0}
\\
\bra d M_{jk}, d \overline{M_{\ell m}} \ket_t
&=\delta_{j \ell} \delta_{k m} \frac{dt}{N}, 
\quad 1 \leq j, k, \ell, m \leq N, \, t \geq 0.
\label{eq:variation1}
\end{align}

Assume that the initial spectrum of $M(0)$ is simple.
We label the $N$ eigenvalues of $M(0)$ following a rule
and express the eigenvalue process of
$(M(t))_{t \geq 0}$ as a labeled process, 
$\bLambda(t) = (\Lambda_j(t))_{j=1}^N$,
$t \geq 0$: 
In Section \ref{sec:holomorphic} below, we will introduce
a lexicographic order of coordinates on 
$\C \simeq \R^2$. Then we label the initial eigenvalues so that
the ordering 
$\Lambda_1(0) < \Lambda_2(0) < \cdots < \Lambda_N(0)$
is satisfied. 
(As a matter of course, such an ordering cannot be conserved
as time is passing.) 
We define
$\sigma:=\inf\{t \geq 0:
\exists j \not=k, 
\Lambda_j(t)=\Lambda_k(t)\}$.
It can be proved that
$\sigma=\infty$ almost surely
\cite[Lemma A.9]{BD20}
\cite[Proposition 3.2]{Yab20}.
Hence the initial ordering gives a unique labeling
for the elements of eigenvalue process
$\bLambda(t) = (\Lambda_j(t))_{j=1}^N$, 
which is valid forever, $t \geq 0$, almost surely. 

The labeled eigenvalue process $(\bLambda(t))_{t \geq 0}$ induces 
\textit{the right and the left eigenvector processes}, 
$\bR_j(t) = (R_{jk}(t))_{k=1}^N$, 
$\L_j(t)= (L_{jk}(t))_{k=1}^N$,
for each $j=1, \dotsm, N$ such that
\begin{align}
M(t) \bR_j(t) &= \Lambda_j(t) \bR_j(t),
\nonumber\\
\L_j^{\rm t}(t) M(t) &= \Lambda_j(t) \L_j^{\rm t}(t),
\quad t \geq 0.
\label{eq:RLevec}
\end{align}
We impose the \textit{bi-orthogonality relation} for eigenvectors, 
\begin{equation} 
(\L_j(t),\overline{\bR_k(t)}) = \delta_{jk}, \quad 1 \leq j,k \leq N,
\, t \geq 0, 
\label{eq:biortho}
\end{equation}
where $(\u, \v)$ denotes the \textit{Hermitian inner product}, 
\begin{equation}
(\u, \v) := \sum_{j=1}^N u_j \overline{v_j} 
= \u^{\rm t} \, \overline{\v}, 
\label{eq:Hip}
\end{equation}
and 
the norm is denoted by $\| \u \|^2:=(\u, \u)$.
The above equations \eqref{eq:RLevec} and \eqref{eq:biortho} 
allow for the \textit{scale transformation}
\begin{equation}
\bR_j(t) \to c_j(t) \bR_j(t), \quad
\L_j(t) \to \frac{1}{c_j(t)} \L_j(t),
\quad 1 \leq j \leq N,
\label{eq:scale}
\end{equation}
for any non-vanishing time-dependent factors
$c_j(t), 1 \leq j \leq N, t \geq 0$. 
In other words, at each time $t \geq 0$, there is ambiguity
in choosing right and left eigenvectors for each eigenvalue $\Lambda_j(t)$,
$1 \leq j \leq N$.

We define the $N \times N$ matrix-valued process
\begin{equation}
S(t) =(S_{jk}(t))_{1 \leq j, k \leq N}
:= (\bR_1(t) \bR_2(t) \cdots \bR_N (t)), \quad t \geq 0, 
\label{eq:def_S}
\end{equation}
that is,
$S_{jk}(t)=R_{kj}(t), 1 \leq j, k \leq N, t \geq 0$.
Then, its inverse-matrix-valued process is 
\[
S^{-1}(t) = (S^{-1}_{jk}(t))_{1 \leq j, k \leq N}
= (\L_1(t) \L_2(t) \cdots \L_N (t))^{\rm t}
=(L_{jk}(t))_{1 \leq j, k \leq N}, \quad t \geq 0,
\]
since the bi-orthogonality \eqref{eq:biortho} implies
\[
(S^{-1}(t) S(t))_{jk}
=\sum_{\ell=1}^N L_{j \ell}(t) R_{k \ell}(t)
=(\L_j(t), \overline{\bR}_k(t))=\delta_{jk}, \quad
1 \leq j, k \leq N.
\]
With the matrix expression 
$\Lambda(t)={\rm diag}(\Lambda_1(t), \dots, \Lambda_N(t))$, $t \geq 0$
of $(\bLambda(t))_{t \geq 0}$, 
\eqref{eq:RLevec} gives
\begin{equation}
M(t) = S(t)\Lambda(t)S^{-1}(t), \quad t \geq 0.
\label{eq:MSLambda}
\end{equation}

Remark that if $(M(t))_{t \geq 0}$ was any Hermitian matrix-valued
process, then $(S(t))_{t \geq 0}$ would be a unitary matrix-valued process:
$S^{\dagger}(t) S(t)=I \iff S^{-1}(t)=S^{\dagger}(t)$, $t \geq 0$, 
where 
$S^{\dagger}(t):=\overline{S}^{\rm t}(t)$ 
denotes the Hermitian conjugate of $S(t)$ and
$I$ is the identity matrix.
This means 
$\L_j(t)=\overline{\bR}_j(t)$, $j=1, \dots, N$,
$t \geq 0$.
Then, \eqref{eq:biortho} implies that
there are no overlaps between any
distinct eigenvector processes in the sense, 
$(\bR_j(t), \bR_k(t))=(\L_j(t), \L_k(t))=0$,
$j \not=k$. 
However, here we have considered the non-Hermitian matrix-valued 
process \eqref{eq:Mt}. Hence, if we define
\begin{align}
A(t) =(A_{jk}(t))_{1 \leq j, k \leq N}:=S^{\dagger}(t) S(t) &\iff 
A_{kj}(t)=(\bR_j(t), \bR_k(t)), \, 1 \leq j, k \leq N, 
\nonumber\\
A^{-1}(t)=(A^{-1}_{jk}(t))_{1 \leq j, k \leq N}
:=S^{-1}(t)(S^{-1})^{\dagger}(t) &\iff
A^{-1}_{jk}(t)=(\L_j(t), \L_k(t)), \, 1 \leq j, k \leq N, 
\label{eq:A}
\end{align}
$t \geq 0$, then
$A(t)$ and $A^{-1}(t)$ are not equal to $I$ in general. 
We see that for  $1 \leq j \not=k \leq N$, 
$A_{jk}(t)$ (resp. $A^{-1}_{jk}(t)$)
represents time-evolution of the
overlap between the $j$-th and the $k$-th 
right (resp. left) eigenvector processes
associated with the distinct eigenvalues $\Lambda_j(t)$ and
$\Lambda_k(t)$, $t \geq 0$.
If we consider products of $A^{-1}_{jk}(t)$ and
$A_{kj}(t)$, $1 \leq j, k \leq N$, 
we can obtain the processes which are invariant under
the scale transformation \eqref{eq:scale}.
\begin{df}
\label{thm:def_of_O}
The eigenvector-overlap process
is defined by the $N \times N$ matrix-valued process as
\begin{equation}
\cO(t)=(\cO_{jk}(t))_{1 \leq j, k \leq N}, \quad t \geq 0
\label{eq:overlap1}
\end{equation}
with
\begin{align}
\cO_{jk}(t) &:= A^{-1}_{jk}(t) A_{kj}(t)
\nonumber\\
&=(\L_j(t), \L_k(t)) (\bR_j(t), \bR_k(t)),
\quad 1 \leq j, k \leq N, \, t \geq 0.
\label{eq:overlap2}
\end{align}
\end{df}
By the definitions of 
the Hermitian inner product \eqref{eq:Hip}
and \eqref{eq:overlap2}, we see that
\[
\cO_{jk}(t)=
\overline{(\L_k(t), \L_j(t))}
\overline{(\bR_k(t), \bR_j(t))}
=\overline{\cO_{kj}(t)}, \quad 
1 \leq j, k \leq N, \,
t \geq 0.
\]
That is, $\cO(t)$ is Hermitian;
$\cO^{\dagger}(t)=\cO(t), t \geq 0$.
Eigenvector-overlaps play important roles in
a variety of fields in mathematics and physics.
See 
\cite{ATTZ20,BNST17,BD20,Burda14,Burda15,BF24,CM98,Fyo18,FM02,FO22,FS12,FT21,GW18,JNNPZ99,MC00,Yab20} 
and references therein.

The statistics of eigenvalues
and the eigenvector-overlaps of the random matrices 
in the Ginibre ensemble has been extensively
studied 
\cite{ATTZ20,BD20,BF24,CM98,For10,Fyo18,FT21,MC00,Meh04}.
Consider the complex Ginibre ensemble
with matrix size $N$
and variance $t/N$.
We can show that the eigenvalue point process
on $\C$ is simple a.s. at any $t>0$, 
and that it converges weakly to the Dirac measure
concentrated on the origin of $\C$ 
as $t \downarrow 0$.
Hence, if 
the present matrix-valued BM starts from
the null matrix $M(0)=O$, then, 
at each time $t \geq 0$
the distribution of $M(t)$ will be identified with
the complex Ginibre ensemble.
Here we study the process starting from
a deterministic but arbitrary matrix $M(0)$.

\subsection{
SDEs with scale-transformation invariance
}
\label{sec:SDEs}

By Lemma \ref{lem:holomorphic} given 
in Section \ref{sec:holomorphic} below, it is verified that 
we have maps from $M(t)$ to $\C$, 
$\Lambda_j(t) = \Phi_j(M(t))$, 
$R_{jk}(t) = \Psi_{jk}(M(t))$, and 
$L_{jk}(t) = \widehat{\Psi}_{jk}(M(t))$, $1 \leq j, k \leq N$, $t \geq 0$, 
such that $\Phi_{j}$, $\Psi_{jk}$, $\widehat{\Psi}_{jk}$, $1 \leq j, k \leq N$
are holomorphic.
Since $(M(t))_{t \geq 0}$ is a matrix-valued BM
\eqref{eq:Mt},
$(\Lambda_j(t))_{t \geq 0}$, $(R_{jk}(t))_{t \geq 0}$,
$(L_{jk}(t))_{t \geq 0}$, $1 \leq j, k \leq N$
are all local martingales in the time period
such that there is no degeneracy in eigenvalues 
$\bLambda(t)=(\Lambda_j(t))_{j=1}^N$ of $M(t)$
(that is, $M(t) \in \C^{N^2} \setminus \Omega$ in the notation
in Section \ref{sec:holomorphic}).

Hence, we can apply It\^o's formula 
to \eqref{eq:MSLambda} to derive 
the stochastic differential equations (SDEs)
for the eigenvalue process $(\bLambda(t))_{t \geq 0}$
\cite{BD20,GW18},
\begin{equation}
d \Lambda_j(t)= (S^{-1}(t) dM(t) S(t))_{jj}, \quad 1 \leq j \leq N, 
\, t \geq 0,
\label{eq:SDE_Lambda1}
\end{equation}
with the cross-variations 
\begin{align}
&\bra d \Lambda_j, d \Lambda_k \ket_t =0,
\quad 
\bra d \overline{\Lambda_j}, d \overline{\Lambda_k} \ket_t =0,
\nonumber\\
&\bra d \Lambda_j, d \overline{\Lambda_k} \ket_t
=\frac{\cO_{jk}(t)}{N} dt,
\quad 1 \leq j, k \leq N, \, t \geq 0.
\label{eq:SDE_Lambda2}
\end{align}
(See Section \ref{sec:A1} for derivation.)
As a matter of course, the equations 
\eqref{eq:SDE_Lambda1} and \eqref{eq:SDE_Lambda2}
are invariant under the scale transformation \eqref{eq:scale},
since they concern the eigenvalue process.

The ambiguity in determination of the eigenvector processes
associated with $(\bLambda(t))_{t \geq 0}$,
which is due to
the covariance of the system
\eqref{eq:RLevec} keeping
\eqref{eq:biortho}
under the scale transformation
\eqref{eq:scale} as mentioned above, 
is revealed in SDEs as follows: 
If we apply It\^o's formula to \eqref{eq:MSLambda},
we obtain
\begin{equation} 
dS_{j k}(t) = S_{j k}(t) dU_{kk}(t) 
+ \sum_{1 \leq \ell \leq N: \ell \neq k} S_{j \ell}(t)
\dfrac{(S^{-1}(t)dM(t)S(t))_{\ell k}}{\Lambda_k(t)-\Lambda_{\ell}(t)},
\quad 1 \leq j, k \leq N, \, t \geq 0,
\label{eq:SDE_S1}
\end{equation}
where $dU(t):=S^{-1}(t) dS(t), t \geq 0$. 
(See Section \ref{sec:A2} for derivation.)
This system of SDEs is not invariant under the scale transformation \eqref{eq:scale}. 
Moreover, we see that the $N$ processes 
$(U_{jj}(t))_{t \geq 0}$, $1 \leq j \leq N$
are unable to be determined 
by the system of SDEs \eqref{eq:SDE_Lambda1}
and \eqref{eq:SDE_S1}. 
In order to eliminate such indeterminate variables,
Grela and Warcho{\l} \cite{GW18} 
impose the following constraint, 
\begin{equation} \label{eq:S-1iidSii}
S^{-1}_{jj}(t) dS_{jj}(t) = 
- \sum_{1 \leq k \leq N: k \neq j} S^{-1}_{jk}(t)dS_{kj}(t),
\quad 1 \leq j \leq N, \,  t \geq 0.
\end{equation}
By the definition of $dU(t), t \geq 0$,
it is obvious that \eqref{eq:S-1iidSii} is equivalent with
$dU_{jj}(t) \equiv 0$, $1 \leq j \leq N, t \geq 0$.
We claim, however, that 
this constraint depends on time and then
we have to change the choice of eigenvectors
with each passing moment 
for keeping a given eigenvalue process 
$(\bLambda(t))_{t \geq 0}$.
Such a procedure will break the analyticity
of the eigenvector processes and hence
application of It\^o's formula could not be justified.

The first main result of this paper is that
the proper system of SDEs 
should be given by the pair of
the SDEs \eqref{eq:SDE_Lambda1} with \eqref{eq:SDE_Lambda2} for the
eigenvalue process
and the following SDEs for 
the eigenvector-overlap process.
For a square matrix $\sm$, ${\rm per} \, \sm$ denotes its permanent.
When the matrix is $2 \times 2$, 
$\displaystyle{\sm=\begin{pmatrix}
m_{11} & m_{12} \\
m_{21} & m_{22} \\
\end{pmatrix}}$, 
then 
${\rm per} \,\sm=m_{11} m_{22}+m_{12} m_{21}$.

\begin{thm}
\label{thm:SDE}
The eigenvector-overlap process $(\cO(t))_{t \geq 0}$ 
given by Definition 1.1
satisfies the following system of SDEs, 
\begin{equation}
d\cO_{jk}(t)
= d \cM^{\cO}_{jk}(t)
+\frac{2}{N} 
\sum_{\substack{1 \leq \ell, m \leq N: \cr \ell \not=j, m \not=k}}
\frac{ \cO_{jk}(t) \cO_{\ell m}(t) +\cO_{\ell k}(t) \cO_{jm}(t)}
{(\Lambda_j(t)-\Lambda_{\ell}(t))
(\overline{\Lambda_k(t)}-\overline{\Lambda_m(t)})} dt,
\label{eq:SDEofdOjk}
\end{equation}
$1 \leq j, k \leq N, t \geq 0$
with
\begin{small}
\begin{align}
&d \cM^{\cO}_{jk}(t)
= \sum_{1 \leq \ell \leq N: \ell \not=j} 
\frac{A_{jk}^{-1}(t)A_{k \ell}(t)
(S^{-1}(t)dM(t) S(t))_{\ell j} 
+ A_{k j}(t)A^{-1}_{\ell k}(t)(S^{-1}(t) dM(t) S(t))_{j \ell}}
{\Lambda_j(t) - \Lambda_{\ell}(t)} \nonumber\\
&\quad + \sum_{1 \leq \ell \leq N: \ell \not=k}
\frac{A^{-1}_{jk}(t) A_{\ell j}(t) 
\overline{(S^{-1}(t) dM(t) S(t))_{\ell k}}
+A_{kj}(t) A^{-1}_{j \ell}(t) 
\overline{(S^{-1}(t) dM(t) S(t))_{k \ell}}}
{\overline{\Lambda_k(t)}-\overline{\Lambda_{\ell}(t)}},
\label{eq:mar_O}
\end{align}
\end{small}
where $(A(t))_{t \geq 0}$ is defined by \eqref{eq:A}.
The increments of quadratic variations 
are given by
\begin{align}
\bra d \cO_{jk}, d \cO_{jk} \ket_t
&=\frac{2 \cO_{jk}(t)}{N}
\sum_{\substack{1 \leq \ell, m \leq N: \cr \ell \not=j, m \not=k}}
\frac{1}{(\Lambda_j(t)-\Lambda_{\ell}(t))
(\overline{\Lambda_k(t)}-\overline{\Lambda_m(t)})}
\nonumber\\
&\times
{\rm per}
\begin{pmatrix}
A_{k \ell}(t) A_{m j}(t) & A_{k \ell}(t)A_{m j}(t)+A_{kj}(t)A_{m \ell}(t) \\
A^{-1}_{\ell k}(t) A^{-1}_{j m}(t) & A^{-1}_{\ell k}(t) A^{-1}_{j m}(t) 
+A^{-1}_{jk}(t)A^{-1}_{\ell m}(t)
\end{pmatrix}
dt,
\label{eq:var_Ojk}
\end{align}
$1 \leq j, k \leq N, t \geq 0$.
\end{thm}
(The proof is given in Section \ref{sec:proof_SDE}.)

It is obvious by Definition 1.1
that at each time $t \geq 0$ the eigenvector-overlap matrix
$\cO(t)=(\cO_{jk}(t))_{1 \leq j, k \leq N}$ is invariant under
the scale transformation \eqref{eq:scale}.
In the SDEs \eqref{eq:SDEofdOjk} 
for $\cO_{jk}(t)$, $1 \leq j, k \leq N$, $t \geq 0$, 
however, the local martingale terms 
\eqref{eq:mar_O} and their quadratic variations
\eqref{eq:var_Ojk} cannot be expressed
only using the eigenvalue process and the
eigenvector-overlap process.
Nevertheless, we can verify that 
they are all invariant under \eqref{eq:scale}.
In particular, 
the invariance of \eqref{eq:var_Ojk}
under \eqref{eq:scale} can be shown
explicitly as explained below.
We define
\begin{align*}
O_{jk\ell m}(t) &:= A^{-1}_{jk}(t) A_{k \ell}(t)
A^{-1}_{\ell m}(t) A_{m j}(t)
\nonumber\\
&= \cO_{\ell m j k}(t),
\quad 1 \leq j, k, \ell, m \leq N,
\quad t \geq 0.
\end{align*}
Then, it is obvious that they are invariant
under \eqref{eq:scale}.
Moreover, by the above definition,
the following 
\textit{decomposition formulas} 
are readily proved, 
\[
\cO_{jkjm}(t)=\cO_{jk}(t) \cO_{jm}(t),
\quad
\cO_{jk \ell k}(t)=\cO_{jk}(t) \cO_{\ell k}(t),
\quad 1 \leq j, k, \ell, m \leq N,
\quad t \geq 0.
\]
We can show that \eqref{eq:var_Ojk} is rewritten as
\[
\bra d \cO_{jk}, d \cO_{jk} \ket_t
=\frac{2 \cO_{jk}(t)}{N}
\sum_{\substack{1 \leq \ell, m \leq N: \cr \ell \not=j, m \not=k}}
\frac{\cO_{jk \ell m}(t)+2 \cO_{jm}(t) \cO_{\ell k}(t)
+\cO_{j m \ell k}(t)}
{(\Lambda_j(t)-\Lambda_{\ell}(t))
(\overline{\Lambda_k(t)}
-\overline{\Lambda_m(t)})} dt,
\]
$1 \leq j, k \leq N$, $t \geq 0$.
We will call $\{\cO_{jk \ell m}(t)\}_{1 \leq j, k, \ell, m \leq N}$,
$t \geq 0$
the \textit{generalized eigenvector-overlap 
processes} associated with $M(t)$, $t \geq 0$.

We will prove the following statement 
in Section \ref{sec:proof_invariance}.

\begin{thm}
\label{thm:invariance}
The SDEs \eqref{eq:SDEofdOjk} with \eqref{eq:mar_O}
for the eigenvector-overlap process $(\cO(t))_{t \geq 0}$ 
is indeed invariant under the scale transformation
\eqref{eq:scale}.
In other words,
although the right and the left 
eigenvector processes
cannot be determined uniquely for a given 
eigenvalue process $(\bLambda(t))_{t \geq 0}$,
the system of SDEs for
the eigenvector-overlap process $(\cO(t))_{t \geq 0}$
is uniquely determined. 
\end{thm} 

\subsection{
Two types of point processes and 
the regularized Fuglede--Kadison determinant}
\label{sec:log_process}

One of the motivations of the present paper
is to make dynamical extensions of
the study by Chalker and Mehlig \cite{CM98,MC00} 
on the eigenvalue and eigenvector statistics
in non-Hermitian random matrix ensembles.
In their study, spacial dependence
of the density of eigenvalues and
the local average of diagonal elements of
overlap matrix were discussed.
Here we define two types of measure-valued processes
so that the expectations of them give
the time-dependent density functions 
$\rho_N$ and $\cO_N$
(see Definition \ref{thm:densities}
in Section \ref{sec:PDE}).

The first one is for the eigenvalue process
$(\bLambda(t))_{t \geq 0}$,
\begin{equation}
\Xi(t, \cdot) := \frac{1}{N} \sum_{j=1}^N \delta_{\Lambda_j(t)}(\cdot),
\quad t \geq 0,
\label{eq:Xi1}
\end{equation}
where $\delta_{\zeta}(\cdot)$ denotes the Dirac measure
concentrated on a point $\zeta \in \C$;
$\delta_{\zeta}(\{z\}) = 1$ if $z =\zeta$,
$\delta_{\zeta}(\{z\}) = 0$ if $z \not= \zeta$.
At each time $t \geq 0$,
$\Xi(t, \cdot)$ gives the normalized empirical measure
of the eigenvalues of $M(t)$.
Another process is defined by
\begin{equation}
\Theta(t, \cdot) := \frac{1}{N^2} \sum_{j=1}^N \cO_{jj}(t) 
\delta_{\Lambda_j(t)}(\cdot),
\quad t \geq 0.
\label{eq:Theta1}
\end{equation}
Although
$\supp \Xi(t, \cdot)=\supp \Theta(t, \cdot)$,
$t \geq 0$, in $(\Theta(t, \cdot))_{t \geq 0}$
each Dirac measure on $\Lambda_j(t)$
is weighted by the $j$-th diagonal element
$\cO_{jj}(t)$, $1 \leq j \leq N$, $t \geq 0$.
Notice that we adopt the scale factor $N^2$ 
instead of $N$ in \eqref{eq:Theta1}. 
(See the footnote 1 on page 426 of \cite{Burda15}.)
By the definition \eqref{eq:overlap2} 
and the Cauchy--Schwarz inequality, 
we can show that \cite{BNST17}, \cite[Sect. 52]{TE05}
\[
\cO_{jj}(t)=\|\L_j(t)\|^2 \| \bR_j(t) \|^2 \ge 1,
\quad 1 \leq j \leq N, \, t \geq 0.
\]

The support of $(\Xi(t))_{t \geq 0}$ defines the 
\textit{time-dependent point process} of eigenvalues on $\C$.
A new kind of time-dependent point process is defined by
$(\Theta(t))_{t \geq 0}$, in which each point
$\Lambda_j(t)$ on $\C$ carries a positive value
$\cO_{jj}(t)$, as a mass or a positive charge
of a point-mass discussed in physics.
As a matter of fact, 
diagonal overlaps $\cO_{jj}, 1 \leq j \leq N$ 
of eigenvectors of non-Hermitian matrices and operators 
have appeared and played important roles in many fields of
mathematics and physics. 
In the numerical analysis, 
the square root of the diagonal overlaps,
$\sqrt{\cO_{jj}}, 1 \leq j \leq N$ 
are known as the 
\textit{eigenvalue condition numbers}.
See \cite[Sects. 35 and 52]{TE05}, 
\cite{BNST17,BD20,Fyo18,FM02,FO22,FS12,MKS25,MKS26,Yab20} 
and references therein.

The time evolution of these empirical measures of
point processes are related with the time-dependent
random fields associated with the 
\textit{Fuglede--Kadison} (FK) \textit{determinant}
of $(M(t))_{t \geq 0}$.
For a regular ({\it i.e.}, invertible) matrix $\sm \in \C^{N^2}$, 
$|\sm|:=(\sm^{\dagger} \sm)^{1/2}$ gives
a positive-definite Hermitian matrix. 
Then the FK determinant 
\cite{FK52} \cite[Chapter 11]{MS17}
is defined by
$\Delta(\sm)=\exp[ \Tr (\log |\sm|)]$.
(See Section \ref{sec:matrices} below for the definitions
and calculus of matrix-valued functions,
\textit{e.g.}, $|\sm|$, $\log |\sm|$.)
If we consider $\sm-z I$ with $z \in \C$, 
instead of $\sm$,
it is not invertible when $z$ coincides with
any eigenvalue $\lambda_j$ of $\sm$, $1 \leq j \leq N$. 
In this case, we introduce a non-zero complex variable
$w \in \C^{\times} := \C \setminus \{0\}$ 
in addition to $z \in \C$
and define the \textit{regularized FK determinant}
\cite{BNST17,BSS18,Burda14,Burda15,MS17} by
\begin{equation}
\Delta_{w}(\sm-z I)
:=\exp \Big[
\Tr \log 
\{(\sm^{\dagger}- \zbar I)(\sm-z I)+ |w|^2 I \}^{1/2}
\Big] \in [0, \infty).
\label{eq:FK2}
\end{equation}
The FK determinant of the matrix 
$\sm-z I$ is then defined by the limit
$\Delta(\sm-z I):=\lim_{w \to 0} \Delta_{w}(\sm-z I) \in [0, \infty)$. 
The following properties are satisfied \cite{MS17}:
\begin{description}
\item{(i)} \,
$\Delta(\sm_1 \sm_2)=\Delta(\sm_1) \Delta(\sm_2),
\quad \forall \sm_1, \sm_2 \in \C^{N^2}$,
\item{(ii)} \,
$\Delta(\sm)=\Delta(\sm^{\dagger})=\Delta(|\sm|), \quad 
\forall \sm \in \C^{N^2}$,
\item{(iii)} \,
$\Delta(I)=1$,
\item{(iv)} \,
$\Delta(c \sm)=|c| \Delta(\sm), \quad
\forall c \in \C, \, \sm \in \C^{N^2}$.
\end{description}

We consider 
the logarithm of \eqref{eq:FK2} and define
\begin{align}
\psi(z, w; \sm) &:= \frac{1}{N} \log \Delta_w(\sm-z I)
\nonumber\\
&= \frac{1}{2N} \Tr \log \widehat{\sh}(z, w; \sm) \in \R,
\quad (z, w) \in \C \times \C^{\times}, \, t \geq 0, 
\label{eq:psi_zwm}
\end{align}
where
\begin{equation}
\widehat{\sh}(z, w; \sm)
:= (\sm^{\dagger}-\zbar I)(\sm- zI) + |w|^2 I.
\label{eq:hat_h}
\end{equation}
It is also written as
\begin{equation}
\psi(z, w; \sm)
=\frac{1}{2N} \log \det \widehat{\sh}(z, w; \sm), 
\quad (z, w) \in \C \times \C^{\times}, \, t \geq 0.
\label{eq:psi2}
\end{equation}
(See Lemma \ref{lem:Meq} given below.)

In the present paper, 
we consider the case such that
$\sm$ in \eqref{eq:FK2} is given by 
the stochastic process $(M(t))_{t \geq 0}$,
and study the random fields 
and random measures defined as
follows.
Notice that these quantities depend
on the matrix size $N \in \N$ of 
$(M(t))_{t \geq 0}$ by definition.
For simplicity of notatoin, however, 
we will suppress any indication of $N$ 
until Definition \ref{thm:densities}
in Section \ref{sec:PDE}, where we will discuss
the $N \to \infty$ limit. 
Let $m$ be the Lebesgue measure on $\C$; 
$m(dz):= dx dy$ for $z=x+i y$ with $x, y \in \R$. 
The Laplacian 
with respect to $z \in \C$ is denoted as
$\displaystyle{\nabla^2_z:= 4 \frac{\partial^2}{\partial z \partial \zbar}}$.
\begin{df}
\label{thm:random_fields_measures}
\begin{description}
\item{\rm (i)} 
The time-dependent random field
associated with the FK determinant
is defined by
\begin{align}
\Delta_w(M(t)-zI)
&=\exp \left[ \Tr \log
\Big\{ (M^{\dagger}(t)-\zbar I) (M(t)-z I) + |w|^2 I \Big\}^{1/2} \right]
\nonumber\\
&=
\sqrt{\det
\Big\{ (M^{\dagger}(t)-\zbar I) (M(t)-z I) + |w|^2 I \Big\}}, 
\quad t \geq 0. 
\label{eq:FK_rf}
\end{align}
\item{\rm (ii)} 
The squared and the logarithmic variations 
of \eqref{eq:FK_rf} are defined 
on the two-dimensional complex space
$(z, w) \in \C \times \C^{\times}$ by
\begin{align}
\Delta^2(z,w; t)
&:= (\Delta_w(M(t)-z I))^2 \quad
\nonumber\\
&= \exp \left[ \Tr \log
\Big\{ (M^{\dagger}(t)-\zbar I) (M(t)-z I) + |w|^2 I \Big\} \right]
\nonumber\\
&=  \det
\Big\{ (M^{\dagger}(t)-\zbar I) (M(t)-z I) + |w|^2 I \Big\}, 
\quad t \geq 0, 
\label{eq:sqFK_rf}
\end{align}
and
\begin{align}
\Psi(z,w; t)
&:=\psi(z,w; M(t))
\nonumber\\
&= \frac{1}{2N} \Tr \log \Big\{
(M^{\dagger}(t)-\zbar I)(M(t)-zI) +|w|^2 I \Big\}
\nonumber\\
&= \frac{1}{2N} \log \det \Big\{
(M^{\dagger}(t)-\zbar I)(M(t)-zI) +|w|^2 I \Big\}, 
\quad t \geq 0,
\label{eq:logFK_rf}
\end{align}
respectively. 
\item{\rm (iii)}
Depending on $w \in \C^{\times}$ and $t \geq 0$, 
the random measures on $\C$ are defined as
\begin{align}
\mu^{\Lambda}_{w}(t, dz)
&:= \frac{2}{\pi} 
\frac{\partial^2 \Psi(z,w; t)}{\partial z \partial \zbar} m(dz)
=\frac{1}{2 \pi} \Big(\nabla_z^2 \Psi(z, w; t) \Big) m(dz),
\label{eq:Brown1}
\\
\mu_{w}^{\cO}(t, dz)
&:= \frac{4}{\pi} \left| 
\frac{\partial \Psi(z, w; t)}{\partial w} \right|^2 m(dz).
\label{eq:Brown2}
\end{align}
\end{description}
\end{df}

Let $\cM_{\rm c}(\C)$ be the set of all Borel 
non-negative measures on $\C$ with compact support equipped with the weak topology.
Let $\cB_{\rm c}(\C)$ be the set of all bounded 
measurable complex functions on $\C$
with compact support. 
For $\mu \in \cM_{\rm c}(\C)$ and 
$\phi \in \cB_{\rm c}(\C)$, 
we set the pairing, 
\begin{equation}
\bra \mu, \phi \ket
:= \int_{\C} \phi(z) \mu(dz).
\label{eq:pair}
\end{equation}

\begin{prop}
\label{prop:psi_Xi_Theta}
For any $\phi \in \cB_{\rm c}(\C)$, 
\begin{align}
\lim_{w \to 0}
\bra \mu^{\Lambda}_w(t, \cdot), \phi(\cdot) \ket
&= \bra \Xi(t, \cdot), \phi(\cdot) \ket, 
\label{eq:psi_Xi}
\\
\lim_{w \to 0}
\bra \mu^{\cO}_w(t, \cdot), \phi(\cdot) \ket
&= \bra \Theta(t, \cdot), \phi(\cdot) \ket, 
\quad t \geq 0.
\label{eq:psi_Theta}
\end{align}
\end{prop}
The statement which is equivalent with
the above proposition was reported in
physics literature \cite{Burda14,Burda15,JNNPZ99}. 
See also \cite{BNST17}. 
We will give a mathematical proof in Section \ref{sec:measures}
with the evaluation given in Section \ref{sec:estimate}. 

The second main result 
of the present paper is the following 
\textit{stochastic partial differential equations} (SPDEs)
for the time-dependent random fields 
on $\C \times \C^{\times}$,
$(\Delta_w(M(t)-zI))_{t \geq 0}$,
$(\Delta^2(z, w; t))_{t \geq 0}$,  and
$(\Psi(z, w; t))_{t \geq 0}$, 
generated by $(M(t))_{t \geq 0}$. 

\begin{thm}
\label{thm:psi_SDE}
The time-dependent random field 
associated with the regularized FK determinant
of $(M(t))_{t \geq 0}$ 
and its squared 
and the logarithmic variations
given by Definition \ref{thm:random_fields_measures} (ii)
satisfy the following SPDEs,
\begin{align}
&d \Delta_w(M(t)-z I)
= d \cM^{\Delta}(z, w; t)
\nonumber\\
& \qquad \qquad +\frac{1}{2N}
\left\{ 
\frac{\partial^2 \Delta_w(M(t)-zI)}{\partial w \partial \wbar}
+ \frac{3}{\Delta_w(M(t)-zI)}
\left| \frac{\partial \Delta_w(M(t)-zI)}{\partial w}
\right|^2 \right\} dt, 
\label{eq:SDE_FK}
\\
&d \Delta^2(z, w; t)
= d \cM^{\Delta^2}(z, w; t)
+ \frac{1}{4N} \nabla^2_w \Delta^2(z, w; t) dt,
\label{eq:SDE_sqFK}
\\
&d \Psi(z, w; t) = d \cM^{\Psi}(z, w; t)
+ 2 \left|
\frac{\partial \Psi(z, w; t)}{\partial w} \right|^2 dt,
\quad (z, w) \in \C \times \C^{\times}, t \geq 0, 
\label{eq:SDE_psi1}
\end{align}
where $\nabla^2_w$ denotes the Laplacian 
with respect to $w \in \C^{\times}$:
$\displaystyle{\nabla^2_w:=4 \frac{\partial^2}{\partial w \partial \wbar}}$.
Here $(\cM^{\Delta^2}(\cdot, \cdot; t))_{t \geq 0}$,
$(\cM^{\Delta}(\cdot, \cdot; t))_{t \geq 0}$, 
and
$(\cM^{\Psi}(\cdot, \cdot; t))_{t \geq 0}$ are 
local martingales satisfying the relations
\begin{align}
d \cM^{\Delta}(z, w; t)
&=N \Delta_w(M(t)-z I) \, d \cM^{\Psi}(z, w; t), 
\nonumber\\
d \cM^{\Delta^2}(z, w; t)
&=2 N \Delta^2(z, w; t) \, d \cM^{\Psi}(z, w; t), 
\label{eq:MMrelation}
\end{align}
with the equation
\begin{align}
d \cM^{\Psi}(z, w; t)
&= \frac{1}{2N} \left\{
\Tr \Big[ \widehat{\sh}^{-1}(z, w; M(t)) (M^{\dagger}(t)-\zbar I)
 d M(t) \Big] \right.
\nonumber\\
& \quad  \quad \, \, \, \, \left.
+\Tr \Big[ (M(t) - zI) \widehat{\sh}^{-1}(z, w; M(t))
d M^{\dagger}(t) \Big] \right\}, 
\label{eq:SDE_psi2}
\end{align}
$(z, w) \in \C \times \C^{\times}, t \geq 0$, 
where $\widehat{\sh}$ is given by \eqref{eq:hat_h}.
The quadratic variation of \eqref{eq:SDE_psi2}
is given by
\begin{equation}
\bra d \cM^{\Psi}(z,w; \cdot), d \cM^{\Psi}(z, w; \cdot)
\ket_t=
\frac{1}{4N^2} \nabla^2_w \Psi(z, w; t) dt,
\quad (z, w) \in \C \times \C^{\times}, \, t \geq 0.
\label{eq:SDE_psi4}
\end{equation}
\end{thm}
\vskip 0.3cm
(The proof is given in Section \ref{sec:psi_SDE}.)

By Definition \ref{thm:random_fields_measures} (iii)
and Proposition \ref{prop:psi_Xi_Theta},
the following is obtained as a corollary of Theorem \ref{thm:psi_SDE}.

\begin{cor}
\label{cor:relation}
For any $\phi \in \cB_{\rm c}(\C)$,
\begin{equation}
d \bra \mu^{\Lambda}_w(t), \phi \ket
=d \cM^{\Lambda}_{w, \phi} (t)
+ \frac{1}{4} \bra \nabla^2_z \mu^{\cO}_w(t), \phi \ket dt,
\quad w \in \C^{\times}, \, t \geq 0,
\label{eq:relation1}
\end{equation}
where $(\cM^{\Lambda}_{w, \phi}(t))_{t \geq 0}$, $w \in \C^{\times}$
is a local martingale given by
\[
\cM^{\Lambda}_{w, \phi}(t)= \frac{1}{2 \pi}
\int_{\C} \phi(z) \Big( \nabla^2_z \cM^{\Psi}(z, w; t) \Big) m(dz). 
\]
In the limit $w \to 0$ of \eqref{eq:relation1}, 
we obtain the SPDE
\begin{equation}
d \bra \Xi(t), \phi \ket
= d \cM^{\Lambda}_{0, \phi}(t)
+\frac{1}{4} \bra \nabla^2_z \Theta(t), \phi \ket dt,
\quad t \geq 0.
\label{eq:relation2}
\end{equation}
\end{cor}
\vskip 0.3cm

\noindent
Here we notice that the SPDE \eqref{eq:relation2}
relates $(\Xi(t))_{t \geq 0}$ and $(\Theta(t))_{t \geq 0}$
and this equation 
can be directly proved from the SDEs
for the eigenvalue process \eqref{eq:SDE_Lambda1}
with \eqref{eq:SDE_Lambda2}
(see Section \ref{sec:proof_relation}). 

\subsection{
PDEs obtained by averaging SPDEs
}
\label{sec:PDE}
We consider the expectation of 
the \textit{logarithmic regularized FK-determinant random-field} 
\eqref{eq:logFK_rf}, which we write as 
\begin{equation}
\bra \Psi \ket (z, w; t) := \E \Big[\Psi(z,w; t) \Big]
=\E \Big[ \psi(z, w; M(t)) \Big],
\quad (z, w) \in \C \times \C^{\times}, \, t \geq 0.
\label{eq:potential0}
\end{equation}
In the context of 
\textit{free probability theory} \cite{BSS18,MS17}, 
the normalized trace 
$\displaystyle{\tau(\sx)=\frac{1}{N} \Tr \sx}$
is called the \textit{tracial state} for $\sx$
which is in a suitable matrix algebra. 
When we consider the filtered probability space 
$(\Omega, \cF, \{\cF_t\}_{t \geq 0}, \P)$
for the matrix-valued process $(X(t))_{t \geq 0}$, 
the tracial state is defined by 
$\displaystyle{\tau(X(t))=\E\left[ \frac{1}{N} \Tr X(t) \right]}$, $ t \geq 0$.
Hence \eqref{eq:potential0} is written as
\begin{equation}
\bra \Psi \ket (z, w; t)
= \tau \Big[ 
\log \{ (M^{\dagger}(t)-\zbar I)(M(t)-z I) + |w|^2 I \}^{1/2}
\Big], 
\label{eq:potential1}
\end{equation}
$(z, w) \in \C \times \C^{\times}, t \geq 0$. 
Then we see that the expectation of \eqref{eq:Brown1} gives
a (time-dependent) deterministic measure 
$\bra \mu^{\Lambda}_w \ket
=\bra \mu^{\Lambda}_w \ket(t, \cdot)$ on $\C$:
\[
\bra \mu^{\Lambda}_w \ket (t, dz)
:=\frac{1}{2 \pi} \Big(\nabla^2_z 
\bra \Psi \ket (z,w;t) \Big) m(dz), \quad t \geq 0.
\]
This is regarded
as the (time-dependent) 
\textit{regularized Brown measure}
which has been extensively studied in free probability theory
\cite{BSS18,BYZ24,DHK22,HH22a,HH22b,HZ19,MS17}.

Next we consider the expectation of the
\textit{squared regularized FK-determinant random-field} \eqref{eq:sqFK_rf},
which we write as
\[
\bra \Delta^2 \ket (z,w; t)
:= \E \big[ \Delta^2(z, w; t) \Big],
\quad (z, w) \in \C \times \C^{\times}, \, t \geq 0.
\]
Since the SDE \eqref{eq:SDE_sqFK} for $\Delta^2$
given in Theorem \ref{thm:psi_SDE} is linear,
we obtain the following diffusion equation for $\bra \Delta^2 \ket$
with respect to $(t, w) \in [0, \infty) \times \C^{\times}$, 
\begin{equation}
\frac{\partial}{\partial t}
\bra \Delta^2 \ket (z, w; t)
=\frac{1}{4N} \nabla_w^2 
\bra \Delta^2 \ket (z, w; t),
\quad (z, w) \in \C \times \C^{\times}, \, t \geq 0.
\label{eq:diffusion}
\end{equation}

\vskip 0.3cm
\noindent \textbf{Remark 1} \,
In \cite{Burda14,Burda15} and the physics literature cited therein,
the deterministic function 
$\bra \Psi \ket (z, w; t)$ given by \eqref{eq:potential1} 
has been introduced as a potential function for the system
following the argument based on ``electrostatic'' analogy. 
Burda et al. \cite{Burda14,Burda15} also studied the
time-dependent deterministic field
$\bra \Delta^2 \ket (z,w; t)$,
which was denoted by $D(z, w;t)$, 
$(z, w) \in \C \times \C^{\times}, t \geq 0$. 
They reported that if we perform the path-integral calculation 
using Grassmann variables, the
diffusion equation \eqref{eq:diffusion} is derived. 
We do not think that it will be easy to justify such calculation
mathematically.
We find, however, if we average 
the SPDE \eqref{eq:SDE_sqFK} in 
Theorem \ref{thm:psi_SDE} by deleting the local martingale term, 
we obtain \eqref{eq:diffusion}. 
We need a mathematical justification to allow
such exchange of the averaging and the
differentiation with respect to the auxiliary 
complex variable $w$, 
which will be a future problem.
\vskip 0.3cm
For the two types of point processes
\eqref{eq:Xi1} and \eqref{eq:Theta1},
we define the two types of \textit{time-dependent 
density functions} as follows.
\begin{df}
\label{thm:densities}
With respect to the Lebesgue measure $m(dz)$, 
$z \in \C$,
the time-dependent density functions
$\rho_N(t,z)$ and $\cO_N(t,z)$
are defined, respectively, 
for the point processes
$(\Xi(t, \cdot))_{t \geq 0}$ and
$(\Theta(t, \cdot))_{t \geq 0}$ by
\begin{align*}
\E[\bra \Xi(t, \cdot), \phi(\cdot) \ket]
&= \int_{\C} \phi(z) \rho_N(t,z) m(dz),
\nonumber\\
\E[\bra \Theta(t, \cdot), \phi(\cdot) \ket]
&= \int_{\C} \phi(z) \cO_N(t,z) m(dz),
\quad \forall \phi \in \cB_{\rm c}(\C),
\quad t \geq 0. 
\end{align*}
\end{df}

By Proposition \ref{prop:psi_Xi_Theta} 
with Definition \ref{thm:random_fields_measures} (iii)
and \eqref{eq:potential0}, 
we have the equalities,
\begin{align}
\rho_N(t, z) &= \frac{1}{2 \pi} \nabla^2_z 
\bra \Psi \ket (z, 0; t), 
\label{eq:density_rho}
\\
\cO_N(t, z) &= \frac{4}{\pi}
\E \left[
\left| \frac{\partial \Psi}{\partial w} (z, 0; t) \right|^2
\right], 
\quad z \in \C, \, t \geq 0, 
\label{eq:density_O}
\end{align}
where 
$\displaystyle{\left. \frac{\partial \Psi}{\partial w}(z, 0; t):=
\frac{\partial \Psi(z,w;t)}{\partial w} \right|_{w=0}}$.
If we take expectations of both sides 
of \eqref{eq:relation2} in Corollary \ref{cor:relation}, 
the local martingale term vanishes and 
we obtain the PDE, 
\begin{equation}
\frac{\partial \rho_N(t,z)}{\partial t}
=\frac{1}{4} \nabla^2_z \cO_N(t,z),
\quad z \in \C, \, t \geq 0.
\label{eq:PDE4}
\end{equation}
We introduce the \textit{current field}
by
\[
j_N(t, z):= - \frac{\partial \cO_N(t,z)}{\partial \zbar},
\quad z \in \C, \, t \geq 0.
\]
Then \eqref{eq:PDE4} is written as
\[
\frac{\partial \rho_N(t, z)}{\partial t}
+\frac{\partial j_N(t, z)}{\partial z}=0, 
\quad z \in \C, \, t \geq 0,
\]
which will be regarded as
the \textit{equation of continuity}.
Hence the density function $\cO_N(t, z)$ 
for the point process \eqref{eq:Theta1}
can be regarded as the potential function
for the current field $j_N(t, z)$
associated with $\rho_N(t, z)$, $t \geq 0, z \in \C$.

\vskip 0.3cm
\noindent \textbf{Remark 2} \,
Let
\begin{equation}
\phi_N(z, w; t) := \frac{1}{2N} 
\log \bra \Delta^2 \ket (z, w; t),
\quad (z, w) \in \C \times \C^{\times}, \, t \geq 0.
\label{eq:phi}
\end{equation}
Since the solution of the diffusion equation \eqref{eq:diffusion}
under the initial condition
$\bra \Delta^2 \ket (z, w; 0)$, $(z ,w) \in \C \times \C^{\times}$ is
given by
\[
\bra \Delta^2 \ket (z, w; t) 
=\frac{N}{\pi t} \int_{\C}
\exp \left( -N \frac{|w-w'|^2}{t} \right) 
\bra \Delta^2 \ket (z, w'; 0) m(d w'),
\]
the $N \to \infty$ limit of \eqref{eq:phi} is evaluated by the
Hopf--Lax formula \cite{Burda15}, 
\begin{align}
\phi(z, w; t)
&:= \lim_{N \to \infty} \phi_N(z, w; t)
\nonumber\\
&= \lim_{n \to \infty}
\frac{1}{n} \log \left[
\frac{n}{2 \pi t} \int_{\C}
\exp \left\{ n \left(
\phi(z, w'; 0) - \frac{|w-w'|^2}{2t} \right) \right\}
m(dw') \right]
\nonumber\\
&=\max_{w'} \left( \phi(z, w'; 0)
-\frac{|w-w'|^2}{2t} \right),
\quad
(z, w) \in \C \times \C^{\times}, \, t \geq 0,
\label{eq:max}
\end{align}
provided that
the limit $\phi(z, w;0)=\lim_{N \to \infty} \phi_N(z, w; 0)$ exists.
Define
\[
v(z, w; t) := \frac{\partial \phi(z,w;t)}{\partial w},
\quad (z, w) \in \C \times \C^{\times}, \, t \geq 0.
\]
If the maximum in \eqref{eq:max} is achieved for
$w'=w_{\ast}$, the extremum condition of \eqref{eq:max}
is 
\[
\overline{v(z, w_{\ast}; 0)}
+\frac{w-w_{\ast}}{2t}=0.
\]
It gives
\[
w_{\ast}= w+2t \overline{v(z, w_{\ast}; 0)}, 
\]
for which $v(z,w; t)=v(z, w_{\ast}; 0)$, $t \geq 0$.
This implies the following functional equation, 
\[
v(z,w;t)=v(z, w+2 t \overline{v(z, w; t)}; 0).
\]
This is the unique solution of the 
\textit{complex Burgers equation in the inviscid limit},
\begin{equation}
\frac{\partial v(z, w; t)}{\partial t}
=2 \frac{\partial |v(z, w; t)|^2}{\partial w},
\quad (z, w) \in \C \times \C^{\times}, \, t \geq 0,
\label{eq:Burgers}
\end{equation}
under the initial function $v(z, w; 0)$,
$(z, w) \in \C \times \C^{\times}$ \cite{Burda15}. 
For the PDEs \eqref{eq:diffusion}, \eqref{eq:PDE4}, and \eqref{eq:Burgers}, 
which have been obtained here by averaging the SPDEs with
deleting the local martingale terms, 
see \cite{BYZ24,DHK22,HH22a,HH22b,HZ19}. 
See also the item (1) in Section \ref{sec:future} below. 
\vskip 0.3cm

The present paper is organized as follows.
In Section \ref{sec:preliminaries}, we only consider 
deterministic matrices, 
where we show 
the formulas that we use in probabilistic 
and stochastic setting in the present paper.
The SDEs for the eigenvalue process \eqref{eq:SDE_Lambda1}
with \eqref{eq:SDE_Lambda2} are derived in
Section \ref{sec:A1}.
In Sections \ref{sec:A2} and \ref{sec:A3},
the SDEs for the right eigenvectors \eqref{eq:SDE_S1}
and for the left eigenvectors \eqref{eq:SDE_S-1}
are derived, respectively.
Section \ref{sec:proofs} is devoted to
the proofs for theorems presented in this paper.
Several future problems are listed out in Section \ref{sec:future}. 

\SSC
{Preliminaries}
\label{sec:preliminaries}
\subsection{Holomorphic functions of a matrix}
\label{sec:holomorphic}

We consider a deterministic complex matrix
\[
\sm = (m_{jk})_{1 \leq j, k \leq N} 
= (m_{jk}^{\rR}+im_{jk}^{\rm I})_{1 \leq j, k \leq N}
\in \C^{N^2}, 
\]
where
$m_{jk}^{\rR}, m_{jk}^{\rm I} \in \R, 1 \leq j, k \leq N$.
We define the closed set $\Omega$ in $\C^{N^2}$ by 
\[
\Omega = \Big\{\sm \in \C^{N^2} : 
\ \text{there is degeneracy in the eigenvalues of $\sm$.} \Big\}.  
\]
For complex numbers $z = x+i y, z'=x'+i y' \in \C$ with
$x, y, x', y' \in \R$, 
an ordering $z < z'$ is introduced as
\[
z < z' \iff \text{ $x < x'$ or [$x = x'$ and $y < y'$]}, 
\]
that is, we consider the lexicographic order on $\C \simeq \R^2$. 
We assume $\sm \in \C^{N^2} \setminus \Omega$. 
Then, the eigenvalues 
$(\lambda_1, \dots, \lambda_N) \in \C^N$ 
can be labeled so that they satisfy the ordering
$\lambda_1 < \lambda_2 < \cdots < \lambda_N$. 
Hence, we can define the functions $\Phi_j : \C^{N^2} \setminus \Omega \to \C$: 
\[
\lambda_j = \Phi_j(\sm), \quad 1 \leq j \leq N.
\]
For each eigenvalue $\lambda_j$, $1 \leq j \leq N$, 
we have its right eigenvector $\br_j = (r_{jk})_{k=1}^N$
and its left eigenvector $\bl_j = (l_{jk})_{k=1}^N$ such that
\[
\sm \br_j = \lambda_j \br_j,
\quad
\bl_j^{\rm t} \sm= \lambda_j \bl_j^{\rm t}.
\]
We remark that for each $\lambda_j$,
$1 \leq j \leq N$, the right and the left 
eigenvectors are 
not uniquely determined and 
$\br_j \in \C^N/\C$ and $\bl_j \in \C^{N}/\C$. 
However, we can prove the following. 
See also \cite[Proposition 3.5]{IKK17}, 
\cite[Proposition 4.6]{PP90}.

\begin{lem} 
\label{lem:holomorphic}
The functions 
$\Phi_j : \C^{N^2} \setminus \Omega \ni \sm \mapsto \lambda_j \in \C$,
$1 \leq j \leq N$ are holomorphic. 
In addition, for each eigenvalue $\lambda_j$, $1 \leq j \leq N$, 
there exist vector-valued holomorphic functions 
${\bf \Psi}_{j}=(\Psi_{jk})_{k=1}^N : \C^{N^2} \setminus \Omega \to \C^N$ and
$\widehat{{\bf \Psi}}_{j}=(\widehat{\Psi}_{jk})_{k=1}^N : 
\C^{N^2} \setminus \Omega \to \C^N$ such that 
${\bf \Psi}_{j}(\sm)$ gives one of the right eigenvectors
and $\widehat{{\bf \Psi}}_{j}(\sm)$
gives one of the left eigenvectors
associated with $\lambda_j$, respectively.
\end{lem}
{\it Proof} \ The first claim is shown by the implicit function theorem of holomorphic functions. 
The second claim is given in 
\cite[Section II.4]{Kato80}. 
\qed

\subsection{Scale-transformation invariance}
\label{sec:scale_transform}

Assume that the matrix $\sm  \in \C^{N^2} \setminus \Omega$.
Then a labeling of the eigenvalues 
$(\lambda_1, \dots, \lambda_N)$ is uniquely
determined as mentioned in Section \ref{sec:holomorphic}.
For each $1 \leq j \leq N$, we choose one of the right eigenvectors 
denoted by $\br_j$. Then the set of corresponding left eigenvectors
$(\bl_j)_{j=1}^N$ is determined so that the bi-orthogonality
\begin{equation}
(\bl_j, \overline{\br_k})=\delta_{jk}, \quad 1 \leq j, k \leq N
\label{eq:bi_orthoZ1}
\end{equation}
is satisfied.

We define the matrix 
$\s=(s_{jk}):= (\br_1 \cdots \br_N)$, that is,
$s_{jk}=r_{kj}$, $1 \leq j, k \leq N$.
By the bi-orthogonality \eqref{eq:bi_orthoZ1}, 
$\s^{-1}=(\bl_{1}, \dots, \bl_N)^{\rm t}$. 
We put $\sa =(a_{jk}) := \s^{\dagger} \s$.
The setting is summarized as follow:
\begin{align*}
& \s^{-1} \s= \s \s^{-1} = I,
\quad
 (\s^{-1})^{\dagger} \s^{\dagger} = \s^{\dagger} (\s^{-1})^{\dagger}= I,
\nonumber\\
& \s^{-1} \sm \s=\blambda
:=\diag(\lambda_1, \dots, \lambda_N),
\quad
\s^{\dagger} \sm^{\dagger} (\s^{-1})^{\dagger}=
\blambdabar
=\diag(\overline{\lambda_1}, \dots, \overline{\lambda_N}), 
\nonumber\\
\sa&:=\s^{\dagger} \s,
\quad \sa^{-1}=\s^{-1}(\s^{-1})^{\dagger}.
\end{align*}
The eigenvector-overlap matrix of the
right and the left eigenvectors
$\so=(o_{jk})_{1 \leq j, k \leq N}$ is then defined by
\begin{equation}
o_{jk}:=a^{-1}_{jk} a_{kj},
\quad 1 \leq j, k \leq N.
\label{eq:overlap3}
\end{equation}

With arbitrary non-vanishing factors $c_j$, $1 \leq j \leq N$,
the scale transformation is defined by
\begin{equation}
\br_j \to \widetilde{\br}_j :=c_j \br_j,
\quad
\bl_j \to \widetilde{\bl}_j := \frac{1}{c_j} \bl_j,
\quad 1 \leq j \leq N.
\label{eq:scale2}
\end{equation}
It is obvious by the definitions that
the eigenvalues, the bi-orthogonality \eqref{eq:bi_orthoZ1} of
the right and the left eigenvectors,
and the eigenvector-overlap matrix $\so$
with the elements \eqref{eq:overlap3} are
invariant under the scale transformation \eqref{eq:scale2}. 
On the other hand, 
for an arbitrary complex matrix $\sw \in \C^{N^2}$, 
$({\bf s}^{-1} \sw {\bf s})_{jk}$ 
is not invariant under \eqref{eq:scale2}. 
However, we can prove the following lemma.

\begin{lem} 
\label{lem:invariance}
Let $\sw \in \C^{N^2}$. 
For $1 \leq j, k, \ell \leq N$, $\ell \not=j$, 
\begin{align}
&a_{jk}^{-1} a_{k \ell} (\s^{-1} \sw \s)_{\ell j}, 
\quad a_{k j} a^{-1}_{\ell k} (\s^{-1} \sw \s)_{j \ell},
\label{eq:inv1}
\\
&a^{-1}_{jk} a_{\ell j}\overline{(\s^{-1} \sw \s)_{\ell k}}, \quad a_{kj} a^{-1}_{j \ell}
\overline{(\s^{-1} \sw \s)_{k \ell}}
\label{eq:inv2}
\end{align}
are invariant under \eqref{eq:scale2}. 
\end{lem}
\noindent
{\it Proof} \
Associated with the scale transformation \eqref{eq:scale2}, 
we set $\widetilde{\s}:=(\widetilde{\br}_1 \cdots \widetilde{\br}_N)$ and
$\widetilde{\sa}=(\widetilde{a}_{jk}):=\widetilde{\s}^{\dagger} \widetilde{\s}$. 
We see 
\begin{align} 
\widetilde{a}_{jk}^{-1} &= (\widetilde{\s}^{-1} (\widetilde{\s}^{\dagger})^{-1})_{jk} 
= \sum_{\ell=1}^N \widetilde{s}^{-1}_{j \ell} (\widetilde{s}^{\dagger})^{-1}_{\ell k} 
= \sum_{\ell=1}^N (\widetilde{\bl}_j)_{\ell} (\overline{\widetilde{\bl}_k})_{\ell} 
\nonumber\\
&= \sum_{\ell=1}^N \left( \frac{1}{c_j} \bl_j \right)_{\ell} 
\left( \overline{\frac{1}{c_k} \bl_k} \right)_{\ell} 
= \frac{1}{c_j \overline{c_k}} \sum_{\ell=1}^N (\bl_j)_{\ell} (\overline{\bl_k})_{\ell} 
= \frac{1}{c_j \overline{c_k}} a^{-1}_{jk}.
\label{eq:tildea-1jk}
\end{align}
By the similar calculation, we have
\begin{equation}
\widetilde{a}_{k \ell} = \overline{c_k} c_{\ell} a_{k \ell}. \label{eq:tildeakp}
\end{equation}
Moreover, we can show that
\begin{equation}
(\widetilde{\s}^{-1} \sw \widetilde{\s})_{\ell j} 
= \sum_{1 \leq p, q \leq N} (\widetilde{\bl}_{\ell})_{p} \sw_{p q} (\widetilde{\br}_j)_{q} 
= \frac{c_j}{c_{\ell}} \sum_{1 \leq p, q \leq N} (\bl_p)_{p} \sw_{p q} (\br_j)_{q} 
= \frac{c_j}{c_{\ell}} (\s^{-1} \sw \s)_{\ell j}. \label{eq:tildes-1Ns}
\end{equation}
Combining \eqref{eq:tildea-1jk}--\eqref{eq:tildes-1Ns}, 
we obtain the equality 
$\widetilde{a}_{jk}^{-1 }\widetilde{a}_{k \ell}
(\widetilde{\s}^{-1} \sw \widetilde{\s})_{\ell j} 
= a_{jk}^{-1} a_{k \ell}(\s^{-1} \sw \s)_{\ell j}$ 
for $1 \leq j, k, \ell \leq N$,
$\ell \not=j$.
By the similar calculation, we can also obtain the equality
$\widetilde{a}_{k j} \widetilde{a}^{-1}_{\ell k} 
(\widetilde{\s}^{-1} \sw \widetilde{\s})_{j \ell}
=a_{k j} a^{-1}_{\ell k} (\s^{-1} \sw \s)_{j \ell}$
for $1 \leq j, k, \ell \leq N$,
$\ell \not=j$.
Hence the invariance of the equalities
\eqref{eq:inv1} under \eqref{eq:scale2} is proved. 
We note that $\sa$ is Hermitian. 
Therefore, the above implies the invariance
of the equalities \eqref{eq:inv2} under \eqref{eq:scale2}.
The proof is hence complete.  \qed

\subsection{Matrix-valued functions 
and functions of matrices}
\label{sec:matrices}

Let $\sh$ be an $N \times N$ Hermitian matrix.
We assume that $\sh$ is positive definite, \textit{i.e.}, 
all eigenvalues are positive. 
We write the eigenvalues of $\sh$ as
$\kappa_j>0, 1 \leq j \leq N$ and put
$\bkappa:= \diag(\kappa_1, \dots, \kappa_N)$. 
Then there is an $N \times N$ unitary matrix $\su$
and we have the decomposition of $\sh$ as
\[
\sh = \su \bkappa \su^{\dagger}.
\]

The powers and the logarithm of matrix $\sh$ are then 
defined as follows:
\begin{align*}
\sh^{\alpha}
&:= \su \, \diag(\kappa_1^{\alpha}, \dots, \kappa_N^{\alpha}) \su^{\dagger},
\quad \alpha \in \R,
\nonumber\\
\log \sh 
&:= \su \, \diag(\log \kappa_1, \dots, \log \kappa_N) \su^{\dagger}. 
\end{align*}
By these definitions, it is easy to verify the equalities
$\sh^{\alpha} \sh^{\beta}=\sh^{\alpha+\beta}$, 
$\alpha, \beta \in \R$
with $\sh^0=I$ and
$\log \sh^{\gamma}=\gamma \log \sh$, $\gamma \in \R$.
The following equality is well known and 
readily proved by the above definition.
\begin{lem}
\label{lem:Meq}
Assume that $\sh$ is a positive-definite Hermitian matrix.
Then 
\[
\Tr \log \sh = \log \det \sh.
\]
\end{lem}

Let $\sm=(m_{jk})_{1 \leq j, k \leq N}$ 
be an arbitrary $N \times N$ complex matrix,
$\sm \in \C^{N^2}$. Put \eqref{eq:hat_h}. 
Then $\widehat{\sh}=\widehat{\sh}(z, w; \sm)$ 
is a positive-definite Hermitian matrix
for $(z, w) \in \C \times \C^{\times}$. 
The following equalities are proved by the above definition and 
straightforward calculation. 

\begin{lem}
\label{lem:Meq2}
For \eqref{eq:psi_zwm}, the following equalities are established: 
\begin{align}
&\frac{\partial \psi}{\partial w}
= \frac{1}{2N} \wbar \Tr \widehat{\sh}^{-1},
\quad 
\frac{\partial \psi}{\partial \wbar}
= \frac{1}{2N} w \Tr \widehat{\sh}^{-1}
=\overline{\frac{\partial \psi}{\partial w}},
\nonumber\\
&\nabla^2_w \psi =
4 \frac{\partial^2 \psi}{\partial w \partial \wbar}
= \frac{2}{N} \left\{ \Tr \widehat{\sh}^{-1} 
- |w|^2 \Tr \widehat{\sh}^{-2} \right\}.
\label{eq:Meq5}
\end{align}
Moreover, the following hold: 
\begin{equation}
\frac{\partial \psi}{\partial m_{jk}}
= \frac{1}{2N} ( \widehat{\sh}^{-1} (\sm^{\dagger}-\zbar I))_{kj},
\quad
\frac{\partial \psi}{\partial \overline{m_{jk}}}
= \frac{1}{2N} ((\sm-z I) \widehat{\sh}^{-1} )_{jk},
\quad 1 \leq j, k \leq N. 
\label{eq:Meq6}
\end{equation}
and
\begin{equation}
\sum_{1 \leq j, k \leq N}
\frac{\partial^2 \psi}{\partial m_{jk} \partial \overline{m_{jk}}}
= \frac{1}{2N} |w|^2 (\Tr \widehat{h}^{-1} )^2
= 2N \left| \frac{\partial \psi}{\partial w} \right|^2.
\label{eq:Meq7}
\end{equation}
\end{lem}

\subsection{Derivatives of 
the logarithmic regularized FK-determinant field}
\label{sec:estimate}

The following is a mathematical justification of the 
argument given by Janik et al. \cite{JNNPZ99}.

The matrix $\widehat{\sh}$ given by \eqref{eq:hat_h} is written as
\begin{align*}
\widehat{\sh} 
&= (\s^{-1})^{\dagger} \s^{\dagger} \widehat{\sh} 
\s \s^{-1}
\nonumber\\
&= (\s^{-1})^{\dagger} \s^{\dagger} 
\{(\sm^{\dagger}-\zbar I)(\sm-z I)+|w|^2 I \}
\s \s^{-1}
\nonumber\\
&= (\s^{-1})^{\dagger} \s^{\dagger} 
(\sm^{\dagger} - \zbar I) 
(\s^{-1})^{\dagger} \s^{\dagger} \s \s^{-1}
(\sm- zI)
\s \s^{-1}
+(\s^{-1})^{\dagger} \s^{\dagger} |w|^2 
\s \s^{-1}
\nonumber\\
&= (\s^{-1})^{\dagger} (\blambdabar - \zbar I) \s^{\dagger} \s
(\blambda-z I) \s^{-1}
+ (\s^{-1})^{\dagger} |w|^2 \s^{\dagger} \s \s^{-1}
\nonumber\\
&= (\s^{-1})^{\dagger} \widetilde{\sh} \s^{-1}
\end{align*}
with
\begin{align}
\widetilde{\sh}
=\widetilde{\sh}(z, w; \sm)
&:=(\blambdabar - \zbar I) \sa (\blambda - z I) + |w|^2 \sa
\nonumber\\
&=\Big(
(\overline{\lambda_j}-\zbar) a_{jk} (\lambda_k-z) +|w|^2 a_{jk}
\Big)_{1 \leq j, k \leq N}.
\label{eq:htilde}
\end{align}
Hence, $\psi$ given by \eqref{eq:psi2} is written as
\[
\psi=\psi(z, w; \sm)=\frac{1}{2N} \log \det \widehat{\sh}
=\frac{1}{2N} \{\log \det (\s^{-1})^{\dagger} + \log \det \s^{-1}\}
+ \frac{1}{2N} \log \det \widetilde{\sh}.
\]
Note that $\s$ does not depend on $z$ and $w$.

With \eqref{eq:htilde}, define
$\widetilde{\sh}_0 :=\widetilde{\sh}(z, 0; \sm)$. 
For $1 \leq n \leq N$, $1 \leq j_1 < \cdots < j_n \leq N$, 
let $\widetilde{\sh}_{\sa}^{(j_1, \dots, j_n)}$
be the matrix obtained from $\widetilde{\sh}_0$ by
replacing its $j_k$-th row 
by the $j_k$-th row of the matrix $\sa$ for
$1 \leq k \leq n$, respectively. 
That is, if $\ell \in \{1, \cdots, N\}
\setminus \{j_1, \dots, j_n\}$, then
the $\ell$-th row of $\widetilde{\sh}_{\sa}^{(j_1, \dots, j_n)}$
is $( (\overline{\lambda_{\ell}}-\zbar) a_{\ell 1}(\lambda_1-z),
\cdots, (\overline{\lambda_{\ell}}-\zbar) a_{\ell N}(\lambda_N-z))$,
while the $j_k$-th row is given by 
$(a_{j_k 1}, \dots, a_{j_k N})$ for $1 \leq k \leq n$. 
Then by the multilinearity of determinant, we have
\begin{align}
\det \widetilde{\sh}
&= \det \widetilde{\sh}_0
+ |w|^2 \sum_{j=1}^N \det \widetilde{\sh}^{(j)}_{\sa}
+|w|^4 \sum_{1 \leq j_1 < j_2 \leq N}
\det \widetilde{\sh}_{\sa}^{(j_1, j_2)} +
\cdots +|w|^{2N} \det \sa.
\label{eq:h_expansion}
\end{align}
Here 
\[
\det \widetilde{\sh}_0
=\det_{1 \leq j, k \leq N} 
\Big[ (\overline{\lambda_j}-\zbar) a_{jk} (\lambda_k-z) \Big]
=\det \sa \times \prod_{j=1}^N |\lambda_j-z|^2,
\]
and all coefficients of $|w|^{2n}, 1 \leq n \leq N$
are polynomials of $z$.
Hence, if $z$ is outside the spectrum of $\sm$;
$z \not= \lambda_j, 1 \leq j \leq N$,
$\psi(z, 0; \sm)=\lim_{w \to 0} \psi(z, w; \sm)$ is
harmonic with respect to $z \in \C$, and then
\begin{align}
\lim_{w \to 0} \nabla_z^2 \psi
&= \lim_{w \to 0} 4 \frac{\partial^2 \psi}{\partial z \partial \zbar}
\nonumber\\
&= \frac{2}{N} \frac{\partial^2}{\partial z \partial \zbar}
\Tr \log \Big[ (\sm^{\dagger}-\zbar I)(\sm-zI) \Big]
\nonumber\\
&= \frac{2}{N} \frac{\partial^2}{\partial z \partial \zbar}
\log \det \Big[ (\blambdabar-\zbar I) \sa (\blambda-z I) \Big]
=0,
\quad z \in \C \setminus \{\lambda_1, \dots, \lambda_N\}.
\label{eq:harmonicA1}
\end{align}
Since \eqref{eq:h_expansion} implies that
$\psi$ is an analytic function of $|w|^2$, we can also conclude that 
\begin{align}
\lim_{w \to 0} \frac{\partial \psi}{\partial w}=
\lim_{w \to 0} \frac{\partial \psi}{\partial \wbar}=0,
\quad z \in \C \setminus \{\lambda_1, \dots, \lambda_N\}.
\label{eq:harmonicA2}
\end{align}

Now we choose arbitrarily one eigenvalue of $\sm$ 
and consider the situation such that the variable $z$ is
in the vicinity of an eigenvalue.
Without loss of generality, we can assume that 
chosen eigenvalue is $\lambda_1$.
With the facts \eqref{eq:harmonicA1} and \eqref{eq:harmonicA2}, 
now we assume that 
\[
|\lambda_1-z| \le |w|
\quad \mbox{for sufficiently small $|w|$}
\]
and we will consider the limit $w \to 0$. 
From \eqref{eq:h_expansion} 
and applying Laplace's expansion for $\widetilde{\sh}^{(1)}_{\sa}$, 
we have
\begin{align*}
\det \widetilde{\sh}
&= \det \widetilde{\sh}_0
+ |w|^2 a_{11}\det \widetilde{\sh}^{(1|1)}_0 + f^{(1)}(z, w),
\end{align*}
where
$\widetilde{\sh}^{(1|j)}_0$ is 
the minor of $\widetilde{\sh}_0$
obtained from $\widetilde{\sh}_0$ by
deleting the first row and the $j$-th column,
$1 \leq j \leq N$, and
\begin{align}
f^{(1)}(z, w):=
&|w|^2\sum_{j: j\neq 1}(-1)^{1+j}a_{1j}\det \widetilde{\sh}^{(1|j)}_0
+|w|^2\sum_{j: j\neq 1}^N  \det \widetilde{\sh}^{(j)}_{\sa}
\nonumber\\
& \quad 
+|w|^4 \sum_{1 \leq j_1 < j_2 \leq N}  
\det \widetilde{\sh}_{\sa}^{(j_1, j_2)} +
\cdots +|w|^{2N}  \det \sa.
\label{eq:deff1}
\end{align}
Note that $f^{(1)}(z, w)$ can be seen 
as a polynomial of  $\lambda_1-z, \overline{\lambda_1}-\zbar$ and $|w|$, 
and the smallest degree of $f^{(1)}(z, w)$ 
with respect to these variables is three. 
Since
\begin{align*}
&\det \widetilde{\sh}^{(1|1)}_0
=\det \Big((\overline{\lambda_j}-\zbar)a_{jk}
(\lambda_k-z)\Big)_{2\le j,k \le N}
=\det \sa^{(1|1)}\prod_{j=2}^N |\lambda_j-z|^2,\\
&o_{11} =a_{11}a^{-1}_{11} =a_{11}\frac{\det \sa^{(1|1)}}{\det \sa},
\end{align*}
we have
\begin{align}
\det \widetilde{\sh}
&=\det \sa \times \prod_{j=2}^N |\lambda_j-z|^2 
\Big\{|\lambda_1-z|^2+o_{11}|w|^2+g^{(1)}(z, w)\Big\},
\label{eq:exactform1}
\end{align}
where 
\begin{align*}
g^{(1)}(z, w)
:=\frac{ f^{(1)}(z, w)}{\det \sa \times \prod_{j=2}^N |\lambda_j-z|^2 }.
\end{align*}
Therefore, the logarithmic derivatives of
\eqref{eq:exactform1} with respect to $z$, $w$, and $\wbar$ 
are obtained as follows:
\begin{align}
\frac{\partial \psi}{\partial z}
&=\frac{1}{2N} \frac{\partial}{\partial z} \log \det \widehat{\sh}
=\frac{1}{2N} \frac{\partial}{\partial z} \log \det \widetilde{\sh}\nonumber\\
&=\frac{1}{2N} \frac{\partial}{\partial z} \sum_{j=2}^N
\log |\lambda_j-z|^2
+\frac{1}{2N} \frac{\partial}{\partial z}
\log \Big\{ |\lambda_1-z|^2 + o_{11} |w|^2 +g^{(1)}(z, w) \Big\} \nonumber\\
&=-\frac{1}{2N} \sum_{j=2}^N
\frac{1}{\lambda_j-z}
+ \frac{1}{2N} \frac{-(\overline{\lambda_1}-\zbar)+g^{(1)}_{z}(z, w)}
{|\lambda_1-z|^2 + o_{11} |w|^2+g^{(1)}(z, w)},\nonumber\\
\frac{\partial \psi}{\partial w}
&=\frac{1}{2N} \frac{\partial}{\partial w} \log \det \widetilde{\sh}
= \frac{1}{2N} \frac{o_{11} \wbar+g^{(1)}_{w}(z, w)}{|\lambda_1-z|^2 + o_{11} |w|^2+g^{(1)}(z, w)},\nonumber\\
\frac{\partial \psi}{\partial \wbar}
&=\frac{1}{2N} \frac{\partial}{\partial \wbar} \log \det \widetilde{\sh}
=\frac{1}{2N} \frac{o_{11} w+g^{(1)}_{\wbar}(z, w)}{|\lambda_1-z|^2 + o_{11} |w|^2+g^{(1)}(z, w)},
\label{eq:deriphi}
\end{align}
where 
$\displaystyle{g^{(1)}_{\zeta}:=\frac{\partial g^{(1)}}{\partial \zeta}}$,
$\zeta=z, \zbar, w, \wbar$. 
The above equalities give
\begin{align}
\frac{\partial^2 \psi}{\partial z \partial \zbar}
&= \frac{1}{2N} \frac{o_{11}|w|^2}
{(|\lambda_1-z|^2+o_{11} |w|^2)^2} +r^{(1)}_1(z, w), 
\label{eq:psi_z_zbarm1}
\\
\left| \frac{\partial \psi}{\partial w} \right|^2
&= \frac{1}{4N^2} \frac{o_{11}^2 |w|^2}
{(|\lambda_1-z|^2+o_{11} |w|^2)^2} + r^{(1)}_2(z, w),
\label{eq:psi_w_2m1}
\\
& \hskip 4cm
z \in \C \setminus \{\lambda_1, \dots, \lambda_N\},
\quad |\lambda_1-z| \leq |w|.
\nonumber
\end{align}
For the remainders $r^{(1)}_{\ell}(z, w), \ell=1,2$, we have the following.
Let $B_r(z)$ denote an open ball centered at $z \in \C$ 
with radius $r >0$.
\begin{lem}
\label{lem:remainders}
For $\phi \in \cB_{\rm c}(\C)$, 
\begin{equation}
\lim_{w \to 0} \int_{B_{|w|}(\lambda_1)} \phi(z)r^{(1)}_{\ell}(z, w)m(dz)=0,
\quad \ell=1,2.
\label{eq:ball}
\end{equation}
\end{lem}
\noindent
{\it Proof} \
We show the proof of \eqref{eq:ball} for $\ell=1$. 
By \eqref{eq:deriphi}, we have
\begin{align*}
r^{(1)}_1(z, w)=
&\frac{1}{2N}\frac{g^{(1)}_{z \zbar}(z, w)}{|\lambda_1-z|^2+o_{11}|w|^2
+g^{(1)}(z, w)}\\
&+\frac{1}{2N}\frac{(\lambda_1-z)g^{(1)}_{z}(z, w)
+(\overline{\lambda_1}-\zbar)g^{(1)}_{\zbar}(z, w)-g^{(1)}_{z}(z, w)g^{(1)}_{\zbar}(z, w)}
{\{|\lambda_1-z|^2+o_{11}|w|^2+g^{(1)}(z, w)\}^2}
\\
&-\frac{1}{2N}
\frac{o_{11}|w|^2}{(|\lambda_1-z|^2+o_{11}|w|^2)^2}
\times\frac{g^{(1)}(z, w)\{o_{11}|w|^2+g^{(1)}(z, w)\}}
{\{|\lambda_1-z|^2+o_{11}|w|^2+g^{(1)}(z, w)\}^2}
\\
&+\frac{1}{2N}
\frac{1}{(|\lambda_1-z|^2+o_{11}|w|^2)^2}
\times\frac{g^{(1)}(z, w) |\lambda_1-z|^4}
{\{|\lambda_1-z|^2+o_{11}|w|^2+g^{(1)}(z, w)\}^2}. 
\end{align*}
In \eqref{eq:deff1}, 
\begin{align*}
\det \widetilde{\sh}^{(1|j)}_0
&=(\lambda_1-z)
\prod_{2 \leq k \leq N, k \neq j}(\lambda_k-z)
\prod_{k=2}^N (\overline{\lambda_k}-\zbar)\det \sa^{(1|j)},
\end{align*}
where $\sa^{(1|j)}$ is the
the minor of $\sa$
obtained from $\sa$ by
deleting the first row and the $j$-th column.
Combining this and  the fact that $\det \widetilde{\sh}^{(j)}_{\sa}$ 
has the factor $(\overline{\lambda_1}-\zbar)$, 
we have 
\begin{align*}
|f^{(1)}(z, w)| = \rO(|w|^3), 
\quad z \in B_{|w|}(\lambda_1), \quad w \to 0.
\end{align*}
Since the first and second derivatives of the coefficient of $|w|^2$ in \eqref{eq:deff1} with respect to $z, \zbar$ do not vanish, we have
\begin{align*}
|f^{(1)}_{z}(z, w)|= \rO(|w|^2), \quad
&|f^{(1)}_{\zbar}(z, w)| = \rO(|w|^2), \quad
|f^{(1)}_{z\zbar}(z, w)| 
= \rO(|w|^2), 
\nonumber\\
&\quad z \in B_{|w|}(\lambda_1), \quad w \to 0.
\end{align*} 
By definition, $g^{(1)}(z, w)$ has no singularity for 
$z \in B_{|w|}(\lambda_1)$.
Hence, for sufficiently small $w$ and $z \in B_{|w|}(\lambda_1)$, 
there exist constants $c_j>0, 1 \leq j \leq 6$, such that
\begin{align*}
&|g^{(1)}(z, w)| \le c_1|w|^3,
\quad |g^{(1)}_{z}(z, w)|\le c_2|w|^2,
\quad |g^{(1)}_{\zbar}(z, w)| \le c_3|w|^2,
\quad |g^{(1)}_{z\zbar}(z, w)| \le c_4|w|^2, \\
&\Big||\lambda_1-z|^2+o_{11}|w|^2+g^{(1)}(z, w)\Big|\ge c_5|w|^2,
\quad 
\Big||\lambda_1-z|^2+o_{11}|w|^2 \Big|\ge c_6|w|^2,
\end{align*}
and we obtain
\begin{align*}
|r^{(1)}_1(z, w)|
&\le\frac{1}{2N}
\Big\{
\frac{c_4|w|^2}{c_5|w|^2}
+\frac{|w|\times c_2|w|^2+|w|\times c_3|w|^2+c_2|w|^2 c_3|w|^2}{(c_5|w|^2)^2}\\
& \quad +\frac{o_{11}|w|^2}{(c_6|w|^2)^2} 
\times \frac{c_1|w|^3\times \{ o_{11} |w|^2+c_1|w|^3\} }{(c_5|w|^2)^2}
+\frac{1}{(c_6 |w|^2)^2} \times
\frac{c_1|w|^3 |w|^4}{(c_5 |w|^2)^2}
\Big\}\\
&\le\frac{1}{2N}\Big( C_1+\frac{C_2}{|w|}\Big),
\end{align*} 
for some constants $C_1, C_2>0$.
Therefore, for any $\phi \in \cB_{\rm c}(\C)$, 
$\| \phi \|_{\infty} := \sup_{z \in \C} |\phi(z)| < \infty$, and 
we have
\begin{align*}
\left|\int_{B_{|w|}(\lambda_1)}\phi(z)r^{(1)}_1(z, w)m(dz) \right|
&\le \| \phi \|_{\infty}\int_{B_{|w|}(\lambda_1)}
\frac{1}{2N}\Big( C_1+\frac{C_2}{|w|}\Big)m(dz)\\
&\le \| \phi \|_{\infty}\frac{1}{2N}\Big( C_1+\frac{C_2}{|w|}\Big)
\times \pi |w|^2 
\to 0, \quad \mbox{as $w \to 0$}.
\end{align*}
The assertion of \eqref{eq:ball} for $\ell=2$ is similarly proved.
\qed
\vskip 0.3cm

\SSC
{SDEs for eigenvalue process 
and eigenvector processes}
\label{sec:SDEsA}
For $N \times N$ matrix-valued processes 
$X(t)=(X_{jk}(t))_{1 \leq j, k \leq N}$, 
$Y(t) = (Y_{jk}(t))_{1 \leq j, k \leq N}$, 
$Z(t)=(Z_{jk}(t))_{1 \leq j, k \leq N}$, 
$\bra dX, ZdY \ket_t$ 
denotes an $N \times N$ matrix-valued process, 
whose $(j, k)$-element 
$( \bra dX, Z dY \ket_t)_{jk}$ 
is given by the finite-variation process 
$\sum_{1 \leq \ell, m \leq N} \bra dX_{j \ell}, Z_{\ell m} dY_{m k} \ket_t$, 
$1 \leq j, k \leq N$, $t \geq 0$. 

\subsection{Derivation of \eqref{eq:SDE_Lambda1} 
and \eqref{eq:SDE_Lambda2}}
\label{sec:A1}
We apply It\^o's formula to \eqref{eq:MSLambda} and obtain
\begin{align*}
dM(t) &= dS(t) \Lambda(t) S^{-1}(t) + S(t) d\Lambda(t) S^{-1}(t) + S(t) \Lambda(t) d(S^{-1})(t) \\
&\quad + \bra dS, d\Lambda \ket_t S^{-1}(t) 
+ S(t) \bra d\Lambda, d(S^{-1}) \ket_t + \bra dS, \Lambda d(S^{-1}) \ket_t. 
\end{align*}
The above gives
\begin{align}
&S^{-1}(t)dM(t)S(t) \notag \\
&= dU(t)\Lambda(t) + d\Lambda(t) + \Lambda(t)dV(t) 
 + \bra dU, d\Lambda \ket_t + \bra d\Lambda, dV \ket_t + \bra dU, \Lambda dV \ket_t, \label{eq:Ito}
\end{align}
where
\begin{equation}
dU(t) := S^{-1}(t)dS(t),
\quad
dV(t) := d(S^{-1})(t)S(t).
\label{eq:def_U_V}
\end{equation}
Note that $V(t) \neq U^{-1}(t)$ in general. 
By $S(t)S^{-1}(t) = I$ for any $t \geq 0$, It\^o's formula gives
$S(t)d(S^{-1})(t) = -dS(t) S^{-1}(t)-\bra dS, d(S^{-1}) \ket_t$.
By multiplying $S^{-1}(t)$ from the left and $S(t)$ from the right, we have 
\begin{equation} \label{eq:dS-1S}
dV(t) = -dU(t)-\bra dU, dV \ket_t,
\end{equation}
and hence,
\begin{equation}
\bra dU, dV \ket_t = - \bra dU, dU \ket_t,
\quad 
\bra d\Lambda, dV \ket_t 
= -\bra d\Lambda, dU \ket_t,
\quad
\bra dU, \Lambda dV \ket_t = -\bra dU, \Lambda dU \ket_t. 
\label{eq:dUdV} 
\end{equation}

First we consider the diagonal part of \eqref{eq:Ito}.
Since $\Lambda(t)$ is diagonal, 
the sum of the first and the third terms 
of the RHS of \eqref{eq:Ito} gives
\begin{align*}
(dU(t)\Lambda(t))_{jj} + (\Lambda(t)dV(t))_{jj} 
&= dU_{jj}(t) \Lambda_j(t) + \Lambda_j(t) dV_{jj}(t) 
= -\Lambda_j(t) (\bra dU, dV \ket_t)_{jj} \notag \\
&= \Lambda_j(t) (\bra dU, dU \ket_t)_{jj} 
= \Lambda_j(t) \sum_{k=1}^N \bra dU_{jk}, dU_{kj} \ket_t, 
\end{align*}
$1 \leq j \leq N$, $t \geq 0$, 
where we used \eqref{eq:dS-1S} and the first equality in \eqref{eq:dUdV}.
By the second equality in \eqref{eq:dUdV}, 
the sum of the fourth and the fifth terms
in the RHS of \eqref{eq:Ito} becomes zero:
\[
(\bra dU, d\Lambda \ket_t)_{jj} + (\bra d\Lambda, dV \ket_t)_{jj}
= \bra d\Lambda_j, dU_{jj} \ket_t -\bra d\Lambda_j, dU_{jj} \ket_t 
= 0, \quad 1 \leq j \leq N, \, t \geq 0.
\]
By the last equality in \eqref{eq:dUdV}, 
the last term in the RHS of \eqref{eq:Ito} gives
\[
(\bra dU, \Lambda dV \ket_t)_{jj} 
= \sum_{k=1}^N \Lambda_k(t) \bra dU_{jk}, dV_{kj} \ket_t 
= -\sum_{k=1}^N \Lambda_k(t) \bra dU_{jk}, dU_{kj} \ket_t,
\quad 1 \leq j \leq N, \, t \geq 0.
\]
Hence the diagonal part of \eqref{eq:Ito} gives the equalities, 
\begin{equation} \label{eq:dlam}
d\Lambda_j(t) = (S^{-1}(t)dM(t)S(t))_{jj} 
- \sum_{k:k \neq j} (\Lambda_j(t)-\Lambda_k(t))\bra dU_{jk}, dU_{kj} \ket_t,
\quad 1 \leq j \leq N, \, t \geq 0.
\end{equation}

Next we consider the off-diagonal part of \eqref{eq:Ito}. 
The sum of the first and the third terms of the RHS of \eqref{eq:Ito} gives
\begin{align}
(dU(t)\Lambda(t))_{jk} + (\Lambda(t)dV(t))_{jk} 
&= dU_{jk}(t) \Lambda_k(t) + \Lambda_j(t)dV_{jk}(t) \notag \\
&= (\Lambda_k(t)-\Lambda_j(t))dU_{jk}(t) 
- \Lambda_j(t) (\bra dU, dV \ket_t)_{jk}, 
\label{eq:S-1dSlam}
\end{align}
$1 \leq j \not= k \leq N$, $t \geq 0$, 
where we used \eqref{eq:dS-1S} in the second equality. 
The sum of the last term in the rightmost side of
\eqref{eq:S-1dSlam} and the
$(j, k)$-element of the last term of the RHS of \eqref{eq:Ito}
is written as
\begin{align*} 
- \Lambda_j(t) (\bra dU, dV \ket_t)_{jk}
+ (\bra dU, \Lambda dV \ket_t)_{jk} 
&= \Lambda_j(t) (\bra dU, dU \ket_t)_{jk} 
- (\bra dU, \Lambda dU \ket_t)_{jk} \notag \\
&= \sum_{\ell: \ell \neq j} 
(\Lambda_j(t)-\Lambda_{\ell}(t)) \bra dU_{j \ell}, dU_{\ell k} \ket_t, 
\end{align*}
$1 \leq j \not= k \leq N$, $t \geq 0$, 
where we used the first and the last equalities in \eqref{eq:dUdV}.
On the other hand, we have
\begin{align*}
(\bra dU, d\Lambda \ket_t)_{jk} + (\bra d\Lambda, dV \ket_t)_{jk} 
&= (\bra dU, d\Lambda \ket_t)_{jk} - (\bra d\Lambda, dU \ket_t)_{jk} \notag \\
&= \bra d\Lambda_k-d\Lambda_j, dU_{jk} \ket_t, 
\end{align*}
$1 \leq j \not= k \leq N$, $t \geq 0$, 
where we used the second equality in \eqref{eq:dUdV}.
Therefore, 
the off-diagonal part of \eqref{eq:Ito} gives
\begin{align}
&(\Lambda_k(t)-\Lambda_j(t))dU_{jk}(t) \notag \\
&\quad = (S^{-1}(t)dM(t)S(t))_{jk} 
- \sum_{\ell: \ell \neq j} (\Lambda_j(t)-\Lambda_{\ell}(t)) 
\bra dU_{j \ell}, dU_{\ell k} \ket_t 
- \bra d\Lambda_k-d\Lambda_j, dU_{jk} \ket_t, 
\label{eq:-lamjlamidU}
\end{align}
$1 \leq j \not= k \leq N$, $t \geq 0$, 
which implies 
\begin{align}
(\Lambda_{\ell}(t)-\Lambda_j(t))(\Lambda_k(t)-\Lambda_{\ell}(t)) 
\bra dU_{j \ell}, dU_{\ell k} \ket_t 
= \bra (S^{-1}dMS)_{j \ell}, (S^{-1}dMS)_{\ell k} \ket_t
\label{eq:lami-lamhdU}
\end{align}
for $1 \leq j, k, \ell \leq N$ with $\ell \neq j, k$, $t \geq 0$. 
Applying \eqref{eq:lami-lamhdU} to \eqref{eq:dlam}, we have
\[
d\Lambda_j(t) 
= (S^{-1}(t)dM(t)S(t))_{jj} 
+ \sum_{k: k\neq j}\dfrac{\bra (S^{-1}dMS)_{jk}, (S^{-1}dMS)_{kj} \ket_t}{\Lambda_j(t)-\Lambda_k(t)}, 
\]
$1 \leq j \leq N$, $t \geq 0$. 
By \eqref{eq:variation0} the second term vanishes,
and \eqref{eq:SDE_Lambda1} is obtained.

We find that
\begin{align*}
\bra (S^{-1}dMS)_{jk}, \overline{(S^{-1}dMS)_{m \ell}} \ket_t 
&= \sum_{\alpha, \beta, \nu, \omega} S^{-1}_{j \alpha}(t) S_{\beta k}(t) 
\overline{S^{-1}_{m \nu}(t)} \; \overline{S_{\omega \ell}(t)} \bra
dM_{\alpha\beta}, d \overline{M_{\nu \omega}} \ket_t \notag \\
&= \sum_{\alpha, \beta, \nu, \omega} S^{-1}_{j \alpha}(t) S_{\beta k}(t) \overline{S^{-1}_{m \nu}(t)} \; \overline{S_{\omega \ell}(t)} \frac{1}{N}\delta_{\alpha\nu}\delta_{\beta\omega} dt \notag \\
&= \frac{1}{N} \left( \sum_{\alpha} S^{-1}_{j \alpha}(t) 
\overline{S^{-1}_{m \alpha}(t)} \right) 
\left( \sum_{\beta} \overline{S_{\beta \ell}(t)} 
S_{\beta k}(t) \right) dt \notag \\
&= \frac{1}{N} A^{-1}_{j m}(t)A_{\ell k}(t) dt, \quad
\end{align*}
$1 \leq j, k, \ell, m \leq N$, $t \geq 0$,
where we used \eqref{eq:variation1} and 
the definition \eqref{eq:A} of $(A(t))_{t \geq 0}$.
By \eqref{eq:overlap2}, we obtain \eqref{eq:SDE_Lambda2}.

\subsection{Derivation of \eqref{eq:SDE_S1}}
\label{sec:A2}

From \eqref{eq:-lamjlamidU} and \eqref{eq:lami-lamhdU}, we have
\begin{align*}
&(\Lambda_k(t)-\Lambda_j(t))dU_{jk}(t) \notag \\
&\quad = (S^{-1}(t)dM(t)S(t))_{jk} 
- \sum_{\ell: \ell \neq j, k} 
\dfrac{1}{\Lambda_{\ell}(t)-\Lambda_k(t)} 
\bra (S^{-1}dMS)_{j \ell}, (S^{-1}dMS)_{\ell k} \ket_t \notag \\
&\quad - (\Lambda_j(t)-\Lambda_k(t))\bra dU_{jk}, dU_{kk} \ket_t 
- \bra d\Lambda_k-d\Lambda_j, dU_{jk} \ket_t,
\end{align*}
$1 \leq j \not=k \leq N$, $t \geq 0$. 
Hence, we have
\begin{align*}
dU_{jk}(t) &= 
\dfrac{(S^{-1}(t)dM(t)S(t))_{jk}}{\Lambda_k(t)-\Lambda_j(t)} 
+ \sum_{\ell: \ell \neq j, k} 
\dfrac{\bra (S^{-1}dMS)_{j \ell}, (S^{-1}dMS)_{\ell k} \ket_t}
{(\Lambda_{\ell}(t)-\Lambda_k(t))(\Lambda_j(t)-\Lambda_k(t))}  \notag \\
&\quad + \bra dU_{jk}, dU_{kk} \ket_t 
+ \dfrac{\bra (S^{-1}dMS)_{kk}-(S^{-1}dMS)_{jj}, dU_{jk} \ket_t}
{\Lambda_j(t)-\Lambda_k(t)}, 
\end{align*}
$1 \leq j \not=k \leq N$, $t \geq 0$, 
where we used \eqref{eq:SDE_Lambda1} 
for the last term in the RHS.

\vskip 0.3cm
\noindent \textbf{Remark 3} \,
By Lemma \ref{lem:holomorphic} in Section
\ref{sec:holomorphic} 
and the definition \eqref{eq:def_S}, 
we can say that 
$S_{jk}(t)=\Psi_{jk}(M(t))$, $1 \leq j, k \leq N$
are holomorphic functions of $M(t)$,
which are continuous in $t \geq 0$ a.s.
Hence by It\^o's formula, 
$dS_{jk}(t)$, $1 \leq j, k \leq N$, $t \geq 0$
are increments of local martingales
induced only by $dM_{jk}(t)$, $1 \leq j, k \leq N$, $t \geq 0$, and the SDEs satisfied by $S_{jk}(t), 1 \le j,k \le N$ 
do not include any terms induced by $d\overline{M(t)}$.
Since $dU(t)=(dU_{jk}(t))_{1 \leq j, k \leq N}$, $t \geq 0$
is defined by \eqref{eq:def_U_V},
$dU_{jk}(t)$, $1 \leq j, k \leq N$, $t \geq 0$ are
increments of local martingales 
also induced only by $dM_{jk}(t)$, 
$1 \leq j, k \leq N$, $t \geq 0$.
\vskip 0.3cm
The facts mentioned in Remark 3 and
\eqref{eq:variation0} imply that
\begin{equation}
\bra dU_{jk}, dU_{\ell m} \ket_t=0,
\quad
\bra (S^{-1} dM S)_{jk}, dU_{\ell m} \ket_t=0,
\quad 1 \leq j, k, \ell, m \leq N, \, t \geq 0, 
\label{eq:zero_cross}
\end{equation}
and hence 
\begin{equation}
dU_{jk}(t)
=\sum_{\ell=1}^N S^{-1}_{j \ell}(t) dS_{\ell k}(t) 
= \dfrac{(S^{-1}(t)dM(t)S(t))_{jk}}{\Lambda_k(t)-\Lambda_j(t)},
\quad 1 \leq j \not= k \leq N, \, t \geq 0.
\label{eq:S-1dSDE1}
\end{equation}

By $S(t)S^{-1}(t)=I$, $t \geq 0$, we have 
\begin{equation}
\sum_{j: j \neq k} S_{m j}(t)S^{-1}_{j \ell}(t) 
= \delta_{m \ell} - S_{m k}(t)S^{-1}_{k \ell}(t),
\quad 1 \leq k, \ell, m \leq N, \, t \geq 0.
\label{eq:Eq1}
\end{equation}
Hence, 
\begin{align}
\sum_{j: j \neq k} S_{m j}(t) \sum_{\ell=1}^N S^{-1}_{j \ell}(t) dS_{\ell k}(t) &=
dS_{m k}(t) - S_{m k}(t) \sum_{\ell=1}^N S^{-1}_{k \ell}(t) dS_{\ell k}(t) \notag \\
&= dS_{m k}(t) - S_{m k}(t) dU_{kk}(t),
\quad 1 \leq k, m \leq N, \, t \geq 0.
\label{eq:S-1dSDE2}
\end{align}
We see that \eqref{eq:S-1dSDE1} and \eqref{eq:S-1dSDE2} give
\eqref{eq:SDE_S1}.

\subsection{SDEs for the left eigenvector process}
\label{sec:A3}
From \eqref{eq:Eq1}, we can also obtain the equality
\[
dS^{-1}_{mk}(t)
=S^{-1}_{mk}(t) d V_{mm}(t)
+ \sum_{j: j \not=m} S^{-1}_{j k}(t) dV_{m j}
=-  S^{-1}_{mk}(t) d U_{mm}(t)
- \sum_{j: j \not=m} S^{-1}_{j k}(t) dU_{m j},
\]
$1 \leq k, m \leq N$, $t \geq 0$, 
where $dV(t):=d(S^{-1})(t) S(t)$, $t \geq 0$, and
we used \eqref{eq:dS-1S}.
Using \eqref{eq:S-1dSDE1}, we obtain the following,
\begin{equation}
dS^{-1}_{jk}(t)
= -S^{-1}_{jk}(t) dU_{jj}(t)
- \sum_{1 \leq \ell \leq N: \ell \not=j}
S^{-1}_{\ell k}(t)
\frac{(S^{-1}(t) dM(t) S(t))_{j \ell}}
{\Lambda_{\ell}(t)-\Lambda_j(t)},
\quad 1 \leq j, k \leq N, \, t \geq 0.
\label{eq:SDE_S-1}
\end{equation}

\SSC
{Proofs}
\label{sec:proofs}
\subsection{Proof of Theorem \ref{thm:SDE}}
\label{sec:proof_SDE}

By the definition of $(A(t))_{t \geq 0}$ given by
the first line of \eqref{eq:A}, we have  
\begin{align}
dA(t) = (dS^{\dagger}(t)) S(t) + S^{\dagger}(t) dS(t) 
+ \bra dS^{\dagger}, dS \ket_t, \quad t \geq 0.
\label{eq:SDEofdA}
\end{align}
For $1 \leq j, k \leq N$, the SDEs \eqref{eq:SDE_S1} give
the cross-variation terms as
\begin{align*}
(\bra dS^{\dagger}, dS \ket_t)_{jk} 
&= \sum_{\ell} \bra d\overline{S_{\ell j}}, dS_{\ell k} \ket_t \notag \\
&= \sum_{\ell} \overline{S_{\ell j}(t)} S_{\ell k}(t) \bra d \overline{U_{jj}}, dU_{kk} \ket_t 
+ \sum_{\ell} \overline{S_{\ell j}(t)} \sum_{m: m \neq k} \frac{S_{\ell m}(t)
\bra d \overline{U_{jj}}, (S^{-1}dM S)_{m k} \ket}{\Lambda_k(t)-\Lambda_m(t)} \notag \\
&\quad + \sum_{\ell} S_{\ell k}(t) \sum_{m: m \neq j } 
\frac{\overline{S_{\ell m}(t)} \bra \overline{(S^{-1}dM S)_{m j}}, dU_{kk} \ket_t}
{\overline{\Lambda_j(t)}-\overline{\Lambda_m(t)}} \notag \\
&\quad + \sum_{\ell} \sum_{p: p \neq j} \sum_{q: q \neq k} 
\frac{\overline{S_{\ell p}(t)} S_{\ell q}(t)
\bra \overline{(S^{-1} dM S)_{p j}}, (S^{-1} dM S)_{q k} \ket_t}
{(\overline{\Lambda_j(t)}-\overline{\Lambda_p(t)})(\Lambda_k(t)-\Lambda_q(t))}.
\end{align*}
By \eqref{eq:variation1} and the definition \eqref{eq:A}, we see that
\[
\bra \overline{(S^{-1} dM S)_{p j}}, (S^{-1} dM S)_{q k} \ket_t
= \frac{1}{N} A^{-1}_{qp}(t) A_{jk}(t) dt,
\quad 1 \leq j, k, p, q \leq N, \, t \geq 0, 
\]
and then, we have
\begin{align}
& (\bra dS^{\dagger}, dS \ket_t)_{jk} 
= A_{jk}(t) \bra d \overline{U_{jj}}, dU_{kk} \ket
+ \sum_{\ell: \ell \neq k}A_{j \ell}(t)
\frac{\bra d \overline{U_{jj}}, (S^{-1}dM S)_{\ell k} \ket}
{\Lambda_k(t)-\Lambda_{\ell}(t)} 
\nonumber\\
&\qquad 
+  \sum_{\ell: \ell \neq j} A_{\ell k}(t)
\frac{\bra \overline{(S^{-1}dM S)_{\ell j}}, dU_{kk} \ket}
{\overline{\Lambda_j(t)}-\overline{\Lambda_{\ell}(t)}} 
+\frac{1}{N} \sum_{\ell: \ell \neq j} \sum_{m: m \neq k}
\frac{ A_{\ell m}(t) A_{m \ell }^{-1}(t) A_{jk}(t)}
{(\overline{\Lambda_j(t)}-\overline{\Lambda_{\ell}(t)})
(\Lambda_k(t)-\Lambda_m(t))} dt,
\label{eq:dS*dS}
\end{align}
$1 \leq j, k \leq N$, $t \geq 0$.
On the other hand, by the definition of $(A^{-1}(t))_{t \geq 0}$
given by the second line of \eqref{eq:A}, we have
\begin{align}
dA^{-1}(t) = (dS^{-1}(t))(S^{-1})^{\dagger}(t) 
+ S^{-1}(t) d (S^{-1})^{\dagger}(t) + \bra dS^{-1}, d(S^{-1})^{\dagger} \ket_t,
\quad t \geq 0. 
\label{eq:SDEofdA-1}
\end{align}
By the calculation similar to the above, 
\eqref{eq:SDE_S-1} gives
\begin{align}
&(\bra dS^{-1}, d(S^{-1})^{\dagger} \ket_t)_{jk} 
= \sum_{\ell} \bra dS^{-1}_{j \ell}, d \overline{S^{-1}_{k \ell}} \ket_t \nonumber\\
&=A^{-1}_{jk}(t) \bra dU_{jj}, d \overline{U_{kk}} \ket_t
+ \sum_{\ell: \ell \neq k}A^{-1}_{j \ell}(t) 
\frac{\bra dU_{jj}, \overline{(S^{-1} dM S)_{k \ell}} \ket_t}
{\overline{\Lambda_{\ell}(t)}-\overline{\Lambda_k(t)}} 
\nonumber\\
& \quad
+ \sum_{\ell: \ell \neq j} A^{-1}_{\ell k}(t) 
\frac{\bra (S^{-1} dM S)_{j \ell}, d \overline{U_{kk}} \ket_t}
{\Lambda_\ell(t)-\Lambda_j(t)} 
+ \frac{1}{N}\sum_{\ell: \ell \neq j} \sum_{m: m \neq k} 
\frac{A^{-1}_{\ell m}(t) A^{-1}_{jk}(t) A_{m \ell}(t)}
{(\Lambda_{j}(t) -\Lambda_{\ell}(t))
(\overline{\Lambda_k(t)} -\overline{\Lambda_m(t)})} dt,
\label{eq:dS-1dS*-1}
\end{align}
$1 \leq j, k \leq N, t \geq 0$.
Now we apply It\^o's formula to each element
$\cO_{jk}(t)$, $1 \leq j, k \leq N$, $t \geq 0$
of the eigenvector-overlap process \eqref{eq:overlap2}:
\begin{equation}
d \cO_{jk}(t) = A_{kj} (t) dA^{-1}_{jk}(t) + A^{-1}_{jk}(t) dA_{kj}(t)
+ \bra dA^{-1}_{jk}, dA_{kj} \ket_t =: d \cM^{\cO}_{jk}(t) + d \cN^{\cO}_{jk}(t), 
\label{eq:Ito_O1}
\end{equation}
where
$d \cM^{\cO}_{jk}(t)$ and $d \cN^{\cO}_{jk}(t)$ denote the local martingale part 
and the finite variation part of $d \cO_{jk}(t)$, respectively,
$1 \leq j, k \leq N$, $t \geq 0$. 
By \eqref{eq:SDEofdA} and \eqref{eq:SDEofdA-1}, we see that
\begin{align*}
d\mathcal{M}_{jk}^{\cO}(t) 
&= A_{kj} (t) \Big[(dS^{-1}(t)(S^{-1})^{\dagger}(t))_{jk} 
+(S^{-1}(t) d(S^{-1})^{\dagger}(t))_{jk} \Big] 
\nonumber\\
& \quad
+ A_{jk}^{-1}(t) \Big[(dS^{\dagger}(t)S(t))_{kj} + (S^{\dagger}(t) dS(t))_{kj} \Big] \notag \\
&= A_{kj}(t) \left[ \sum_{\ell} dS^{-1}_{j \ell}(t) \overline{S^{-1}_{k\ell}(t)} 
+ \sum_{\ell} S^{-1}_{j \ell}(t) d \overline{S^{-1}_{k\ell}(t)} \right] 
\nonumber\\
& \quad
+ A^{-1}_{jk}(t) \left[ \sum_{\ell} d \overline{S_{\ell k}(t)} S_{\ell j}(t)
+ \sum_{\ell} \overline{S_{\ell k}(t)} dS_{\ell j}(t) \right],
\quad 1 \leq j, k \leq N, \, t \geq 0.
\end{align*}
If we apply the SDEs \eqref{eq:SDE_S1} and \eqref{eq:SDE_S-1}, 
the above is written as
\begin{equation}
d\mathcal{M}_{jk}^{\cO}(t) 
=d\mathcal{M}_{jk}^{\cO(1)}(t) 
+d\mathcal{M}_{jk}^{\cO(2)}(t)
\label{eq:Y3}
\end{equation}
with 
\begin{align*}
d\mathcal{M}_{jk}^{\cO(1)}(t)  &= A_{kj}(t) 
\left[ \sum_{\ell} (-S^{-1}_{j \ell}(t) dU_{jj}(t)) 
\overline{S^{-1}_{k \ell}(t)} 
+ \sum_{\ell} S^{-1}_{j \ell}(t) (-\overline{S^{-1}_{k \ell}(t)} \; 
d \overline{U_{kk}(t)}) \right] \\
&\quad + A^{-1}_{jk}(t) 
\left[ \sum_{\ell} \overline{S_{\ell k}(t)} \; 
d \overline{U_{kk}(t)} S_{\ell j}(t) 
+ \sum_{\ell} \overline{S_{\ell k}(t)} S_{\ell j}(t) dU_{jj}(t) \right], \\
d\mathcal{M}_{jk}^{\cO(2)}(t) &= A_{kj} (t) 
\left[ \sum_{\ell} \left( -\sum_{m: m \neq j} 
\frac{S^{-1}_{m \ell}(t) (S^{-1}(t) dM(t) S(t))_{j m}}
{\Lambda_{m}(t)-\Lambda_{j}(t)} \right) \overline{S^{-1}_{k \ell}(t)} 
\right. \nonumber\\
& \qquad \qquad \qquad \left.
+ \sum_{\ell} S^{-1}_{j \ell}(t)
\left( -\sum_{m: m \neq k} 
\frac{\overline{S^{-1}_{m \ell}(t)} \;\overline{(S^{-1}(t) dM(t) S(t) )_{k m}}}{\overline{\Lambda_{m}(t)}-\overline{\Lambda_{k}(t)}} \right) \right] \\
&\quad + A^{-1}_{jk}(t) 
\left[ 
\sum_{\ell} 
\sum_{m: m \neq k} 
\frac{\overline{S_{\ell m}(t)} \;\overline{(S^{-1}(t) dM(t) S(t))_{m k}}}
{\overline{\Lambda_{k}(t)}-\overline{\Lambda_{m}(t)}} S_{\ell j}(t)
\right.
\nonumber\\
& \qquad \qquad \qquad \left.
+ \sum_{\ell} \overline{S_{\ell k}(t)} \sum_{m: m \neq j} 
\frac{S_{\ell m}(t) (S^{-1}(t) dM(t) S(t))_{m j}}
{\Lambda_{j}(t)-\Lambda_{m}(t)} \right], 
\end{align*}
$1 \leq j, k \leq N$, $t \geq 0$. 
By the definition \eqref{eq:A}, we see that
\begin{align}
d\mathcal{M}_{jk}^{\cO(1)}(t)
&= A_{kj}(t) (-A^{-1}_{jk}(t) dU_{jj}(t) 
- A^{-1}_{jk}(t) d \overline{U_{kk}(t)})
\nonumber\\
& \quad
+ A^{-1}_{jk}(t) (A_{kj}(t) d \overline{U_{kk}(t)} 
+ A_{kj}(t) dU_{jj}(t)) = 0,
\quad 1 \leq j, k \leq N, \, t \geq 0. 
\label{eq:IY3}
\end{align}
On the other hand, we see
\begin{align}
d\mathcal{M}_{jk}^{\cO(2)}(t)
&= A_{kj}(t) \left[  -\sum_{\ell: \ell \neq j} 
\frac{A^{-1}_{\ell k}(S^{-1}(t) dM(t) S(t))_{j \ell}}{\Lambda_{\ell}(t)-\Lambda_j(t)} 
-\sum_{\ell: \ell \neq k} 
\frac{A^{-1}_{j \ell}(t) \overline{(S^{-1}(t) dM(t) S(t))_{k \ell}}}
{\overline{\Lambda_{\ell}(t)}-\overline{\Lambda_k(t)}} \right] \notag \\
&\quad +A^{-1}_{jk}(t) \left[ \sum_{\ell: \ell \neq k} 
\frac{A_{\ell j}(t) \overline{(S^{-1}(t) dM(t) S(t))_{\ell k}}}
{\overline{\Lambda_k(t)}-\overline{\Lambda_{\ell}(t)}} 
+ \sum_{\ell: \ell \neq j} 
\frac{A_{k \ell}(S^{-1}(t) dM(t) S(t))_{\ell j}}
{\Lambda_j(t)-\Lambda_{\ell}(t)} \right] \notag \\
&= \sum_{\ell: \ell \neq j} \frac{A^{-1}_{jk}(t) A_{k \ell}(t) (S^{-1}(t) dM(t) S(t))_{\ell j} 
+ A_{k j}(t) A^{-1}_{\ell k}(t) (S^{-1}(t) dM(t) S(t))_{j \ell}}
{\Lambda_j(t) - \Lambda_{\ell}(t)} \notag \\
&+ \sum_{\ell: \ell \neq k} 
\frac{A^{-1}_{jk}(t) A_{\ell j}(t) 
\overline{(S^{-1}(t) dM(t) S(t))_{\ell k}}
+A_{kj}(t) A^{-1}_{j \ell}(t) 
\overline{(S^{-1}(t) dM(t) S(t))_{k \ell}}}
{\overline{\Lambda_k(t)}-\overline{\Lambda_{\ell}(t)}}, 
\label{eq:IIY3}
\end{align}
$1 \leq j, k \leq N$, $t \geq 0$. 
By \eqref{eq:Ito_O1}--\eqref{eq:IIY3}, \eqref{eq:mar_O} is obtained.

Next, we calculate $d \cN_{jk}^{\cO}(t)$,
$1 \leq j, k \leq N$, $t \geq 0$. 
From \eqref{eq:SDEofdA} and \eqref{eq:SDEofdA-1}, we have
\begin{align}
d \cN_{jk}^{\cO}(t)
=A_{kj}(t) (\bra dS^{-1}, d(S^{-1})^{\dagger} \ket_t)_{jk}
+A^{-1}_{jk}(t) (\bra dS^{\dagger}, dS \ket_t)_{kj}
+\bra dA^{-1}_{jk}, dA_{kj} \ket_t,
\label{eq:dN}
\end{align}
$1 \leq j, k \leq N$, $t \geq 0$. 
Here we introduce an abbreviation, 
\[dX_{jk}(t):=\frac{(S^{-1}(t) dM(t) S(t))_{jk}}{\Lambda_k(t)-\Lambda_j(t)},
\quad 1 \leq j \not= k \leq N, \, t \geq 0. 
\]
Then, by \eqref{eq:variation0},
\begin{equation}
\bra dX_{jk}, dX_{\ell m} \ket_t =0,
\label{eq:crossvardX0}
\end{equation}
and
by \eqref{eq:variation1}, 
\begin{align}
\bra dX_{jk}, \overline{dX_{\ell m}} \ket_t 
= \frac{1}{N} \frac{A^{-1}_{j \ell}(t) A_{m k}(t)}
{(\Lambda_k(t)-\Lambda_j(t))(\overline{\Lambda_{m}(t)}-\overline{\Lambda_{\ell}(t)})}dt,
\label{eq:crossvardX}
\end{align}
for $1 \leq j, k, \ell, m \leq N$, 
$j \not=k$, $\ell \not= m$, $t \geq 0$. 
Using \eqref{eq:dS-1dS*-1}, 
the first term 
in the RHS of \eqref{eq:dN} is written as
\begin{align}
A_{kj}(t) (\bra dS^{-1}, d(S^{-1})^{\dagger} \ket_t)_{jk} 
&=A_{kj}(t) A^{-1}_{jk}(t) \bra dU_{jj}, d \overline{U_{kk}} \ket_t
+ A_{kj}(t) \sum_{\ell: \ell \neq k} A^{-1}_{j \ell}(t)
\bra dU_{jj}, d \overline{X_{k \ell}} \ket_t \nonumber \\ 
& \, 
+ A_{kj}(t) \sum_{\ell: \ell \neq j} A^{-1}_{\ell k}(t) 
\bra dX_{j \ell}, d \overline{U_{kk}} \ket_t
\nonumber\\
& \,
+ A_{kj}(t) \frac{1}{N}\sum_{\ell: \ell \neq j} \sum_{m: m \neq k} 
\frac{A^{-1}_{\ell m}(t) A^{-1}_{jk}(t) A_{m \ell}(t)}
{(\Lambda_{\ell}(t)-\Lambda_j(t))(\overline{\Lambda_m(t)}-\overline{\Lambda_k(t)})} dt,
\label{eq:fvAdA-1}
\end{align}
$1 \leq j, k \leq N$, $t \geq 0$.
Similarly, using \eqref{eq:dS*dS}, the second term 
in the RHS of \eqref{eq:dN} is written as
\begin{align}
A^{-1}_{jk}(t) (\bra dS^{\dagger}, dS \ket_t)_{kj}
&=A^{-1}_{jk}(t) A_{kj}(t) \bra d \overline{U_{kk}}, dU_{jj} \ket_t 
+ A^{-1}_{jk}(t) \sum_{\ell: \ell \neq j} A_{k \ell}(t) \bra d \overline{U_{kk}}, dX_{\ell j} \ket_t\nonumber\\
& \,
+  A^{-1}_{jk}(t) \sum_{\ell: \ell \neq k} A_{\ell j}(t) 
\bra d \overline{X_{\ell k}}, dU_{jj} \ket_t 
\nonumber\\
& \,
+ A^{-1}_{jk}(t) \frac{1}{N} \sum_{\ell: \ell \neq k} 
\sum_{m: m \neq j} \frac{ A_{\ell m}A_{m \ell}^{-1}(t) A_{kj}(t)}
{(\overline{\Lambda_k(t)}-\overline{\Lambda_{\ell}(t)})
(\Lambda_j(t)-\Lambda_m(t))} dt, 
\label{eq:fvA-1dA}
\end{align}
$1 \leq j, k \leq N$, $t \geq 0$.
By \eqref{eq:A}, 
the last term 
in the RHS of \eqref{eq:dN} is written as
\begin{align*}
&\bra dA^{-1}_{jk}, dA_{kj} \ket_t = \bra (dS^{-1}(S^{-1})^{\dagger})_{jk} 
+ (S^{-1}d(S^{-1})^{\dagger})_{jk}, (dS^{\dagger}S)_{kj} 
+ (S^{\dagger} dS)_{kj} \ket_t \nonumber\\
&= \bra (dS^{-1}(S^{-1})^{\dagger})_{jk}, (dS^{\dagger} S)_{kj} \ket_t 
+ \bra (S^{-1} d(S^{-1})^{\dagger})_{jk}, (S^{\dagger} dS)_{kj} \ket_t \nonumber\\
&= \sum_{p, q} (S^{-1})^{\dagger}_{pk}(t) S_{qj}(t) 
\bra dS^{-1}_{jp}, \overline{dS_{qk}} \ket_t 
+ \sum_{p,q} S^{-1}_{jp}(t) S^{\dagger}_{kq}(t)
\bra \overline{dS^{-1}_{kp}}, dS_{qj} \ket_t
\nonumber\\
&= \sum_{p, q} (S^{-1})^{\dagger}_{pk}(t) S_{qj}(t) 
\left\bra -S^{-1}_{jp}dU_{jj}-\sum_{\alpha: \alpha \neq j}
S^{-1}_{\alpha p} dX_{j \alpha}, \ 
\overline{S_{qk}} \; d \overline{U_{kk}}
+\sum_{\beta: \beta \neq k} 
\overline{S_{q\beta}} \; d \overline{X_{\beta k}} \right\ket_t \nonumber\\
&+ \sum_{p,q} S^{-1}_{jp}(t)  S^{\dagger}_{kq}(t)
\left\bra -\overline{S^{-1}_{kp}} \; d \overline{U_{kk}}
-\sum_{\alpha: \alpha \neq k} \; \overline{S^{-1}_{\alpha p}} 
d \overline{X_{k\alpha}},\  
S_{qj}dU_{jj}
+\sum_{\beta: \beta \neq j}S_{q\beta}dX_{\beta j} \right\ket_t,
\end{align*}
$1 \leq j, k \leq N$, $t \geq 0$, 
where we have applied
\eqref{eq:SDE_S1} and \eqref{eq:SDE_S-1}
using the facts \eqref{eq:zero_cross} and
\eqref{eq:crossvardX0}. 
This can be expanded as
\begin{align*}
& -\sum_{p, q} (S^{-1})^{\dagger}_{pk}(t) S_{qj}(t) S^{-1}_{jp}(t) 
\overline{S_{qk}(t)} \bra dU_{jj},  d \overline{U_{kk}} \ket_t
\nonumber\\
&
-\sum_{p, q} (S^{-1})^{\dagger}_{pk}(t) S_{qj}(t) S^{-1}_{jp}(t) 
\sum_{\beta: \beta \neq k} 
\overline{S_{q\beta}(t)} \bra dU_{jj}, d \overline{X_{\beta k}}\ket_t 
\nonumber\\
&-\sum_{p, q} (S^{-1})^{\dagger}_{pk}(t) S_{qj}(t) 
\overline{S_{qk}(t)}
\sum_{\alpha: \alpha \neq j} S^{-1}_{\alpha p}(t) 
\bra dX_{j \alpha}, d \overline{U_{kk}} \ket_t 
\nonumber\\
&-\sum_{p, q} (S^{-1})^{\dagger}_{pk}(t) S_{qj}(t) 
\sum_{\alpha: \alpha \neq j} S^{-1}_{\alpha p}(t) 
\sum_{\beta: \beta \neq k} \overline{S_{q\beta}(t)} 
\bra dX_{j \alpha}, d \overline{X_{\beta k}} \ket_t
 \nonumber\\
&- \sum_{p,q} S^{-1}_{jp}(t) S^{\dagger}_{kq}(t) 
\overline{S^{-1}_{kp}(t)} S_{qj}(t) 
\bra d \overline{U_{kk}}, dU_{jj} \ket_t
\nonumber\\
&- \sum_{p,q} S^{-1}_{jp}(t) S^{\dagger}_{kq}(t) 
\overline{S^{-1}_{kp}(t)} 
\sum_{\beta: \beta \neq j} S_{q \beta}(t) 
\bra d \overline{U_{kk}}, dX_{\beta j} \ket_t 
\nonumber\\ 
&- \sum_{p,q} S^{-1}_{jp}(t) S^{\dagger}_{kq}(t) S_{qj}(t) 
\sum_{\alpha: \alpha \neq k} \overline{S^{-1}_{\alpha p}(t)} 
\bra d \overline{X_{k\alpha}}, dU_{jj} \ket_t
\nonumber\\
&- \sum_{p,q} S^{-1}_{jp}(t) S^{\dagger}_{kq}(t) 
\sum_{\alpha: \alpha \neq k} 
\overline{S^{-1}_{\alpha p}(t)} 
\sum_{\beta: \beta \neq j} S_{q\beta}(t) 
\bra d \overline{X_{k\alpha}}, dX_{\beta j} \ket_t,
\end{align*}
$1 \leq j, k \leq N$, $t \geq 0$.
By \eqref{eq:A}, the above is written as
\begin{align}
&-A^{-1}_{jk}(t) A_{kj}(t) \bra dU_{jj},  d \overline{U_{kk}} \ket_t 
-A^{-1}_{jk}(t) \sum_{\beta: \beta \neq k} A_{\beta j}(t) 
\bra dU_{jj}, d \overline{X_{\beta k}}\ket_t 
\nonumber\\
&-A_{kj}(t) \sum_{\alpha: \alpha \neq j} 
A^{-1}_{\alpha k}(t) 
\bra dX_{j \alpha}, d \overline{U_{kk}} \ket_t
-\sum_{\alpha: \alpha \neq j} \sum_{\beta: \beta \neq k} 
A^{-1}_{\alpha k}(t) A_{\beta j}(t) 
\bra dX_{j \alpha}, d \overline{X_{\beta k}} \ket_t 
\nonumber\\
&- A^{-1}_{jk}(t) A_{kj}(t) 
\bra d \overline{U_{kk}}, dU_{jj} \ket_t
-A^{-1}_{jk}(t) 
\sum_{\beta: \beta \neq j} A_{k \beta}(t) 
\bra d \overline{U_{kk}}, dX_{\beta j} \ket_t 
\nonumber\\
&- A_{kj}(t) 
\sum_{\alpha: \alpha \neq k} A^{-1}_{j \alpha}(t) 
\bra d \overline{X_{k\alpha}}, dU_{jj} \ket_t
-\sum_{\alpha: \alpha \neq k} 
\sum_{\beta: \beta \neq j} A^{-1}_{j \alpha}(t) 
A_{k \beta}(t) 
\bra d \overline{X_{k\alpha}}, dX_{\beta j} \ket_t, 
\label{eq:fvdA-1dA}
\end{align}
$1 \leq j, k \leq N$, $t \geq 0$.
Here, by \eqref{eq:crossvardX}, 
the fourth and the last terms in 
\eqref{eq:fvdA-1dA} are equal to
\begin{align}
&-\frac{1}{N} \sum_{\alpha: \alpha \neq j} 
\sum_{\beta: \beta \neq k}  
\frac{A^{-1}_{\alpha k} (t)A_{\beta j}(t) A^{-1}_{j\beta}(t) A_{k \alpha}(t)}
{(\Lambda_{\alpha}(t)-\Lambda_j(t))
(\overline{\Lambda_{k}(t)}-\overline{\Lambda_{\beta}(t)})}dt
\nonumber\\
\mbox{and} \quad
&-\frac{1}{N} \sum_{\alpha: \alpha \neq k} 
\sum_{\beta: \beta \neq j}
\frac{A^{-1}_{j \alpha}(t) A_{k \beta}(t) A^{-1}_{\beta k}(t) A_{\alpha j}(t)}
{(\Lambda_j(t)-\Lambda_{\beta}(t))
(\overline{\Lambda_{\alpha}(t)}-\overline{\Lambda_k(t)})}dt, 
\label{eq:sumA}
\end{align}
respectively, $1 \leq j, k \leq N$, $t \geq 0$.
It is easy to see that these two terms
are equal to each other 
by exchanging the indices $\alpha \leftrightarrow \beta$.
Inserting \eqref{eq:fvAdA-1}--\eqref{eq:fvdA-1dA}
with \eqref{eq:sumA}
into \eqref{eq:dN}, we find 
the following cancellation of terms in
$d \cN_{jk}^{\cO}(t)$, $1 \leq j, k \leq N$, $t \geq 0$. 
Here notice that the quadratic variation is symmetric
with respect to the arguments: 
$\bra X, Y \ket_t=\bra Y, X \ket_t$, $t \geq 0$.
For $1 \leq j, k \leq N$, $t \geq 0$, 
\begin{align*}
&\mbox{(the terms including
$\bra dU_{jj}, d \overline{U_{kk}} \ket_t$)}
\nonumber\\
& \qquad 
=A_{kj}(t) A^{-1}_{jk}(t) \bra dU_{jj}, d \overline{U_{kk}} \ket_t
+ A^{-1}_{jk}(t) A_{kj}(t) \bra d \overline{U_{kk}}, dU_{jj} \ket_t 
\nonumber\\
& \qquad 
-A^{-1}_{jk}(t) A_{kj}(t) \bra dU_{jj},  d \overline{U_{kk}} \ket_t
- A^{-1}_{jk}(t) A_{kj}(t) \bra d \overline{U_{kk}}, dU_{jj} \ket_t=0,\\
&\mbox{(the terms including
$\bra dU_{jj}, d \overline{X_{ab}} \ket_t$)} 
\nonumber\\
& \qquad = A_{kj}(t) \sum_{\ell: \ell \neq k} A^{-1}_{j \ell}(t) 
\bra dU_{jj}, d \overline{X_{k \ell}} \ket_t 
+  A^{-1}_{jk}(t) \sum_{\ell: \ell \neq k} A_{\ell j}(t) 
\bra d \overline{X_{\ell k}}, dU_{jj} \ket_t \\
& \qquad -A^{-1}_{jk}(t) \sum_{\beta: \beta \neq k}A_{\beta j}(t)
\bra dU_{jj}, d \overline{X_{\beta k}} \ket_t
- A_{kj}(t) \sum_{\alpha: \alpha \neq k}
 A^{-1}_{j \alpha} \bra d \overline{X_{k\alpha}}, dU_{jj} \ket_t =0,
\\
&\mbox{(the terms including 
$\bra d \overline{U_{kk}}, dX_{ab} \ket_t$)}
\nonumber\\
& \qquad =A_{kj}(t) \sum_{\ell: \ell \neq j} A^{-1}_{\ell k}(t) 
\bra dX_{j \ell}, d \overline{U_{kk}}\ket_t
+ A^{-1}_{jk}(t) \sum_{\ell: \ell \neq j}A_{k \ell}(t)
 \bra d \overline{U_{kk}}, dX_{\ell j} \ket_t \\
& \qquad -A_{kj}(t) \sum_{\alpha: \alpha \neq j} A^{-1}_{\alpha k}(t)
\bra dX_{j \alpha}, d \overline{U_{kk}} \ket_t
-A^{-1}_{jk}(t) \sum_{\beta: \beta \neq j} A_{k \beta}(t)
\bra d \overline{U_{kk}}, dX_{\beta j} \ket_t 
=0.
\end{align*}
And we see that
the remaining terms in $d\cN_{jk}^{\cO}(t)$, $1 \leq j, k \leq N$, $t \geq 0$
are summarized as 
\begin{align*}
&\frac{2}{N} \sum_{\ell: \ell \neq j} \sum_{m: m \neq k} 
\frac{ A_{kj}(t) A^{-1}_{\ell m}(t) A^{-1}_{jk}(t) A_{m \ell}(t)}
{(\Lambda_j(t)-\Lambda_{\ell}(t))(\overline{\Lambda_k(t)}-\overline{\Lambda_m(t)})} dt \\
&
-\frac{2}{N} \sum_{\alpha: \alpha \neq j} \sum_{\beta: \beta \neq k}  
\frac{A^{-1}_{\alpha k}(t) A_{\beta j}(t) A^{-1}_{j\beta}(t) A_{k \alpha}(t)}
{(\Lambda_{\alpha}(t)-\Lambda_j(t))
(\overline{\Lambda_{k}(t)}-\overline{\Lambda_{\beta}(t)})}dt \\
&=\frac{2}{N}\sum_{\ell: \ell \neq j} \sum_{m: m \neq k}
\frac{ \cO_{jk}(t) \cO_{\ell m}(t) +\cO_{\ell k}(t) \cO_{j m}(t)}
{(\Lambda_j(t)-\Lambda_{\ell}(t))(\overline{\Lambda_k(t)}-\overline{\Lambda_m(t)})} dt,
\quad 1 \leq j, k \leq N, \, t \geq 0. 
\end{align*}
That is, the dependence on $dU_{jj}(t)$, $1 \leq j \leq N$, $t \geq 0$, 
disappears completely.
The finite-variation terms
in \eqref{eq:SDEofdOjk} are hence proved. 

The cross-variations of $\cO_{jk}(t)$, $1 \leq j, k \leq N$,
$t \geq 0$ are
calculated using \eqref{eq:variation0} and \eqref{eq:variation1}
and written by the definitions \eqref{eq:A} and \eqref{eq:overlap1} 
with \eqref{eq:overlap2} as
\begin{align*}
&\bra d \cO_{jk}, d \cO_{jk} \ket_t
=\frac{2}{N}\sum_{\ell: \ell \neq j} 
\sum_{m:  m \neq k}\frac{1}{(\Lambda_j(t)-\Lambda_\ell(t))
(\overline{\Lambda_k}(t)-\overline{\Lambda_m(t)})}\\
&\, \times \Big[
A^{-1}_{jk}(t) A_{k \ell}(t) A^{-1}_{jk}(t) A_{m j}(t) A^{-1}_{\ell m}(t) A_{kj}(t)
+A^{-1}_{jk}(t) A_{k \ell}(t) A^{-1}_{j m}(t) A_{kj}(t) A^{-1}_{\ell k}(t) A_{m j}(t)
\\
&\quad+A^{-1}_{\ell k}(t) A_{kj}(t) A^{-1}_{jk}(t) A_{m j}(t) A^{-1}_{j m}(t) A_{k \ell}(t)
+A^{-1}_{\ell k}(t) A_{kj}(t) A^{-1}_{j m}(t) A_{kj}(t) A^{-1}_{jk}(t) A_{m \ell}(t)
\Big] dt 
\\
&=\frac{2 \cO_{jk}(t)}{N}\sum_{\ell: \ell \neq j} 
\sum_{m: m \neq k}\frac{1}{(\Lambda_j(t)-\Lambda_{\ell}(t))
(\overline{\Lambda_k(t)}-\overline{\Lambda_m(t)})}\\
&\qquad \times \Big[
A_{k \ell}(t) A^{-1}_{jk}(t) A_{m j}(t) A^{-1}_{\ell m}(t)
+A_{k \ell}(t) A^{-1}_{j m}(t) A^{-1}_{\ell k}(t) A_{m j}(t) 
\nonumber\\
& \qquad \quad 
+A^{-1}_{\ell k}(t) A_{m j}(t) A^{-1}_{j m}(t) A_{k \ell}(t) 
+A^{-1}_{\ell k}(t) A_{kj}(t) A^{-1}_{j m}(t) A_{m \ell}(t) \Big]dt,
\quad 1 \leq j, k \leq N, \, t \geq 0.
\end{align*}
Using permanents, they are expressed as 
\eqref{eq:var_Ojk}.

\subsection{Proof of Theorem \ref{thm:invariance}}
\label{sec:proof_invariance}

By Lemma \ref{lem:invariance} 
in Section \ref{sec:scale_transform},
we can conclude that the RHS of the SDEs
\eqref{eq:SDEofdOjk} with \eqref{eq:mar_O} 
are invariant under the scale transformation 
\eqref{eq:scale}. 
Notice that the initial value of the eigenvector-overlap process
$\cO(0)=(\cO_{jk}(0))_{1 \leq j, k \leq N}$ 
are independent of the choice of eigenvectors of $M(0)$.
As implied by Lemma \ref{lem:holomorphic}, 
given the non-Hermitian matrix-valued BM
$(M(t))_{t \geq 0}$ 
and its eigenvalue process $(\bLambda(t))_{t \geq 0}$, 
the vector-valued holomorphic functions,
${\bf \Psi}_{j} : \C^{N^2} \setminus \Omega \ni M(t)
\mapsto \bR_j(t) \in \C^N$ and
$\widetilde{\bf \Psi}_{j} : \C^{N^2} \setminus \Omega \ni M(t)
\mapsto \L_j(t) \in \C^N$, $1 \leq j \leq N$, 
$t \geq 0$ cannot be
determined uniquely. 
Nevertheless, 
the eigenvector-overlap process $(\cO(t))_{t \geq 0}$ is uniquely determined.
The proof of Theorem \ref{thm:invariance} is hence complete.

\subsection{Proof of Proposition \ref{prop:psi_Xi_Theta}}
\label{sec:measures}

Assume $\phi \in \cB_{\rm c}(\C)$.
Then by 
\eqref{eq:Brown1} of 
Definition \ref{thm:random_fields_measures} (iii), 
\eqref{eq:pair},
the fact \eqref{eq:harmonicA1}, 
and \eqref{eq:psi_z_zbarm1} with Lemma \ref{lem:remainders}
in Section \ref{sec:estimate}, 
we obtain
\begin{align}
\lim_{w \to 0} \bra \mu^{\Lambda}_w(t, \cdot), \phi(\cdot) \ket
&= \lim_{w \to 0} \int_{\C} \phi(z) 
\frac{1}{N} \sum_{j=1}^N \frac{\cO_{jj}(t) |w|^2}
{\pi (|\Lambda_j(t)-z|^2 + \cO_{jj}(t) |w|^2)^2} m(dz)
\nonumber\\
&= \lim_{\varepsilon \to 0} \int_{\C} \phi(z) 
\frac{1}{N} \sum_{j=1}^N \frac{\varepsilon}
{\pi (|\Lambda_j(t)-z|^2 + \varepsilon)^2} m(dz),
\quad t \geq 0.
\label{eq:Limit1}
\end{align}
For $\phi \in \cB_{\rm c}(\C)$, 
the Cauchy integral formula 
(see, for instance, \cite{Bel16}) gives 
\begin{align*}
\phi(\zeta) &= \lim_{\varepsilon \to 0} 
\frac{1}{2 \pi i} \int_{\C} 
\frac{\partial \phi(z)}{\partial \zbar} \frac{1}{z-\zeta + \varepsilon}
dz \wedge d \zbar
\nonumber\\
&= - \lim_{\varepsilon \to 0} 
\int_{\C} \frac{\partial \phi(z)}{\partial \zbar}
\frac{1}{\pi(z-\zeta+\varepsilon)} m(dz)
\nonumber\\
&=-\lim_{\varepsilon \to 0} 
\int_{\C} \frac{\partial \phi(z)}{\partial \zbar}
\frac{\zbar - \zetabar}{\pi(|z-\zeta|^2+\varepsilon)} m(dz), 
\quad \zeta \in \C.
\end{align*}
By partial integration, it is written as
\[
\phi(\zeta) = \lim_{\varepsilon \to 0}
\int_{\C} \phi(z) 
\left(\frac{\partial}{\partial \zbar}
\frac{\zbar-\zetabar}{\pi(|z-\zeta|^2+\varepsilon)} \right)
m(dz).
\]
Since
\[
\frac{\partial}{\partial \zbar}
\frac{\zbar-\zetabar}{\pi(|z-\zeta|^2+\varepsilon)}
=
\frac{(|z-\zeta|^2+\varepsilon)
-(\zbar-\zetabar)(z-\zeta)}
{\pi(|z-\zeta|^2+\varepsilon)^2}
=\frac{\varepsilon}
{\pi(|z-\zeta|^2+\varepsilon)^2},
\]
we obtain the equality
\[
\phi(\zeta)
= \lim_{\varepsilon \to 0}
\int_{\C} \phi(z)
\frac{\varepsilon}{\pi (|\zeta-z|^2+\varepsilon)^2} m(dz),  
\quad \phi \in \cB_{\rm c}(\C), \, \zeta \in \C.
\]
It gives an expression for the Dirac measure
$\delta_{\zeta}(\cdot), \zeta \in \C$,
\begin{equation}
\bra \delta_{\zeta}, \phi \ket
=\lim_{\varepsilon \to 0} 
\left\bra 
\frac{\varepsilon}{\pi(|\zeta-\cdot|^2+\varepsilon)^2} m(\cdot), 
\phi(\cdot)
\right\ket, \quad
\phi \in \cB_{\rm c}(\C).
\label{eq:Dirac2}
\end{equation}
If we apply this formula to \eqref{eq:Limit1}
and use the definition \eqref{eq:Xi1}, 
then the equality \eqref{eq:psi_Xi} is obtained.

Similarly, 
by the definition \eqref{eq:Brown2}, 
the fact \eqref{eq:harmonicA2}, 
and \eqref{eq:psi_w_2m1} with Lemma \ref{lem:remainders}
in Section \ref{sec:estimate}, we have 
\begin{align*}
\lim_{w \to 0} \bra \mu^{\cO}_w(t, \cdot), \phi(\cdot) \ket
&= \lim_{w \to 0} \int_{\C} \phi(z) 
\frac{1}{N^2} \sum_{j=1}^N \frac{\cO_{jj}(t)^2 |w|^2}
{\pi (|\Lambda_j(t)-z|^2 + \cO_{jj}(t) |w|^2)^2} m(dz)
\nonumber\\
&= \lim_{\varepsilon \to 0} \int_{\C} \phi(z) 
\frac{1}{N^2} \sum_{j=1}^N \cO_{jj}(t)
 \frac{\varepsilon}
{\pi (|\Lambda_j(t)-z|^2 + \varepsilon)^2} m(dz)
\nonumber\\
&= \int_{\C} \phi(z) 
\frac{1}{N^2} \sum_{j=1}^N \cO_{jj}(t) \delta_{\Lambda_j(t)}(dz),
\quad t \geq 0,
\end{align*}
where \eqref{eq:Dirac2} was used.
Hence \eqref{eq:psi_Theta} is proved by
the definition \eqref{eq:Theta1}. 
The proof of Proposition \ref{prop:psi_Xi_Theta} is complete.

\subsection{Proof of Theorem \ref{thm:psi_SDE}}
\label{sec:psi_SDE}

First we prove \eqref{eq:SDE_psi1} with
\eqref{eq:SDE_psi2} and \eqref{eq:SDE_psi4}. 
By It\^o's formula, 
\begin{align*}
d \Psi(z, w; t)
&= \sum_{j=1}^N \sum_{k=1}^N 
\frac{\partial \psi(z, w; \sm)}{\partial m_{jk}} \Big|_{\sm=M(t)} dM_{jk}(t)
+ \sum_{j=1}^N \sum_{k=1}^N 
\frac{\partial \psi(z, w; \sm)}{\partial \overline{m_{jk}}} \Big|_{\sm=M(t)} 
d \overline{M_{jk}(t)}
\nonumber\\
&\quad + \frac{1}{2} \sum_{j=1}^N \sum_{k=1}^N
\frac{\partial^2 \psi(z, w; \sm)}{\partial m_{jk} \partial \overline{m_{jk}} }
\Big|_{\sm=M(t)}
\bra d M_{jk}, d \overline{M_{jk}} \ket_t 
\nonumber\\
&\quad +\frac{1}{2} \sum_{j=1}^N \sum_{k=1}^N
\frac{\partial^2 \psi(z, w; \sm)}{\partial \overline{m_{jk}} \partial m_{jk}}
\Big|_{\sm=M(t)}
\bra d \overline{M_{jk}}, d M_{jk} \ket_t 
\nonumber\\
&= \sum_{j=1}^N \sum_{k=1}^N 
\frac{\partial \psi(z, w; \sm)}{\partial m_{jk}} \Big|_{\sm=M(t)} dM_{jk}(t)
+ \sum_{j=1}^N \sum_{k=1}^N 
\frac{\partial \psi(z, w; \sm)}{\partial \overline{m_{jk}}} \Big|_{\sm=M(t)} 
d \overline{M_{jk}(t)}
\nonumber\\
& \quad + \frac{1}{N} \sum_{j=1}^N \sum_{k=1}^N
\frac{\partial^2 \psi(z, w; \sm)}{\partial m_{jk} \partial \overline{m_{jk}}}
\Big|_{\sm=M(t)} dt,
\quad (z, w) \in \C \times \C^{\times},
\, t \geq 0, 
\end{align*}
where \eqref{eq:variation0} and \eqref{eq:variation1} are used.

Using \eqref{eq:Meq6} of Lemma \ref{lem:Meq2} 
in Section \ref{sec:matrices},
the local martingale term \eqref{eq:SDE_psi2} 
is obtained. 
The finite-variation term in \eqref{eq:SDE_psi1}
is obtained using \eqref{eq:Meq7} of Lemma \ref{lem:Meq2}.

The quadratic variation of the local martingale term
is calculated using \eqref{eq:variation0} and \eqref{eq:variation1} as
\[
\bra d \cM^{\Psi}(z, w; \cdot), d \cM^{\Psi}(z, w; \cdot) \ket_t
= \frac{1}{2N^3} 
\Tr \Big[ (M^{\dagger}(t)-\zbar I)(M(t) - zI) \widehat{h}^{-2}(z, w; M(t))
\Big] dt, 
\]
$(z, w) \in \C \times \C^{\times}$, $t \geq 0$. 
By the definition \eqref{eq:hat_h}, we have
$(M^{\dagger}(t)-\zbar I)(M(t) - zI) 
=\widehat{h}(z, w; M(t))-|w|^2 I$,
and hence the above is written as
\[
\frac{1}{2N^3} 
\Big( \Tr \widehat{h}^{-1}(z, w; M(t))
- |w|^2 \Tr \widehat{h}^{-2}(z, w; M(t)) \Big) dt,
\quad (z, w) \in \C \times \C^{\times},
\, t \geq 0. 
\]
Then \eqref{eq:Meq5} of Lemma \ref{lem:Meq2}
proves \eqref{eq:SDE_psi4}.

Next we apply It\^o's formula to
\begin{align*}
\Delta_w(M(t)-zI) 
&= \exp \Big[ N \Psi(z,w; t) \Big],
\nonumber\\
\Delta^2(z, w; t)
&=(\Delta_w(M(t)- zI))^2
=\exp \Big[ 2N \Psi(z,w; t) \Big],
\quad (z, w) \in \C \times \C^{\times}, \, t \geq 0, 
\end{align*}
using \eqref{eq:SDE_psi1} and \eqref{eq:SDE_psi4}.
Then \eqref{eq:SDE_FK} and \eqref{eq:SDE_sqFK} 
with \eqref{eq:MMrelation}
are obtained. 
The proof of Theorem \ref{thm:psi_SDE} is 
hence complete.

\subsection{Direct proof of \eqref{eq:relation2}
by SDEs \eqref{eq:SDE_Lambda1} with \eqref{eq:SDE_Lambda2}
}
\label{sec:proof_relation}

Let $\phi \in \cB_{\rm c}(\C)$.
We apply It\^o's formula to 
$\displaystyle{
\bra \Xi(t), \phi \ket
=\frac{1}{N} \sum_{j=1}^N \phi(\Lambda_j(t))}$
and obtain
\begin{align}
d \bra \Xi(t), \phi \ket
&=\frac{1}{N} \sum_{j=1}^N
\left\{
\frac{\partial \phi}{\partial z}(\Lambda_j(t))
d \Lambda_j(t)
+ \frac{\partial \phi}{\partial \zbar}
(\Lambda_j(t)) d \overline{\Lambda_j(t)}
+ \frac{1}{2} \frac{\partial^2 \phi}{\partial z^2}(\Lambda_j(t))
\bra d \Lambda_j, d \Lambda_j \ket_t \right.
\nonumber\\
& \qquad \qquad \left.
+ \frac{\partial^2 \phi}{\partial z \partial \zbar}(\Lambda_j(t))
\bra d \Lambda_j, d \overline{\Lambda_j} \ket_t  
+\frac{1}{2} 
\frac{\partial^2 \phi}{\partial \zbar^2}(\Lambda_j(t))
\bra d \overline{\Lambda_j}, d \overline{\Lambda_j} \ket_t
\right\} dt
\nonumber\\
&=\frac{1}{N} \sum_{j=1}^N
\left\{
\frac{\partial \phi}{\partial z}(\Lambda_j(t))
d \Lambda_j(t)
+ \frac{\partial \phi}{\partial \zbar}
(\Lambda_j(t)) d \overline{\Lambda_j(t)}
+ \frac{\partial^2 \phi}{\partial z \partial \zbar}(\Lambda_j(t))
\frac{\cO_{jj}(t)}{N} 
\right\} dt
\nonumber\\
&= d \cM^{\Lambda}_{\phi}(t)
+ \frac{1}{4 N^2} \sum_{j=1}^N \nabla^2_z \phi (\Lambda_j(t))
\cO_{jj}(t) dt, \quad t \geq 0
\label{eq:another}
\end{align}
with a real local martingale, 
\begin{align*}
d \cM^{\Lambda}_{\phi}(t)
&=\frac{1}{N} \sum_{j=1}^N
\left\{
\frac{\partial \phi}{\partial z}(\Lambda_j(t))
(S^{-1}(t) dM(t) S(t))_{jj} \right.
\nonumber\\
& \left. \qquad \qquad 
+ \frac{\partial \phi}{\partial \zbar} (\Lambda_j(t)) 
\overline{(S^{-1}(t) dM(t) S(t))_{jj}}
\right\}, 
\quad t \geq 0.
\end{align*}
Here SDEs \eqref{eq:SDE_Lambda1} with 
\eqref{eq:SDE_Lambda2} were used.
The last term in \eqref{eq:another} is written as
\[
\left\bra \Theta(t), \frac{1}{4} \nabla^2_z \phi 
\right\ket
= \frac{1}{4} \int_{\C}
(\nabla^2_z \phi(z) )
\frac{1}{N^2} \sum_{j=1}^N \cO_{jj}(t)
\delta_{\Lambda_j(t)} (dz),
\quad t \geq 0. 
\]
We perform the integration by parts twice.
Then \eqref{eq:relation2} is obtained. 

\SSC
{Future Problems}
\label{sec:future}

At the end of this paper, we list out
several future problems related to the present study.

\begin{description}
\item{(1)} \,
Assume that for any $N \in \N$,
the initial empirical measures
of the two types of point processes,
$\Xi(0, \cdot)$ and $\Theta(0, \cdot)$, have bounded supports,
and in $N \to \infty$ they converge weakly to the measures
$\nu^{\Lambda}(0, \cdot) \in \cM_{\rm c}(\C)$
and $\nu^{\cO}(0, \cdot) \in \cM_{\rm c}(\C)$, respectively.
For $T>0$, let
$\cC([0,T] \to \cM_{\rm c}(\C))$ denote
the space of continuous processes defined in the
time period $[0, T]$ realized in $\cM_{\rm c}(\C)$.
Then for any arbitrary but fixed $T < \infty$,
we expect the following weak convergence in $N \to \infty$: 
\begin{align*}
& (\Xi(t, \cdot))_{t \in [0, T]} \Longrightarrow (\nu^{\Lambda}(t, \cdot))_{t \in [0, T]}, 
\nonumber\\
& (\Theta(t, \cdot))_{t \in [0, T]} \Longrightarrow (\nu^{\cO}(t, \cdot))_{t \in [0, T]}
\quad \mbox{a.s. in} \, \,  \cC([0, T] \to \cM_{\rm c}(\C)),
\end{align*}
where $\nu^{\Lambda}(t, \cdot)$ and $\nu^{\cO}(t, \cdot)$, $t \geq 0$ 
are time-dependent deterministic measures.
The one-point correlation functions of these limit measures
will be given by
\[
\nu^{\Lambda}(t, dz) = \rho(t, z) m(dz),
\quad
\nu^{\cO}(t, dz)=\cO(t, z) m(dz),
\quad z \in \C, \, t \geq 0.
\]
By definition, $\rho$ and $\cO$ are the $N \to \infty$ limits
of \eqref{eq:density_rho} and \eqref{eq:density_O}, respectively. 
Burda et al. \cite{Burda15} conjectured that
$\cO(t, z)$ shall be obtained by 
\[
\cO(t, z)=\frac{4}{\pi} | v(z, 0; t)|^2,
\quad z \in \C, \, t \geq 0,
\]
where $v(z, 0; t)$ is the $w \to 0$ limit
of the solution $v(z, w; t)$,
$(z, w) \in \C \times \C^{\times}, t \geq 0$ of the
inviscid complex Burgers equation \eqref{eq:Burgers}.
Then $\rho(t, z)$ will be determined as the solution of
the PDE
\[
\frac{\partial \rho(t,z)}{\partial t}
=\frac{1}{4} \nabla^2_z \cO(t,z),
\quad z \in \C, \, t \geq 0,
\]
which is implied by \eqref{eq:PDE4}.
There $\cO(t, z)$ plays the role of potential function
for the current field associated with
$\rho(t, z)$, $t \geq 0, z \in \C$.
Proving the above conjectures 
concerning the $N \to \infty$ limit is a future problem.
See \cite{BYZ24,DHK22,HH22a,HH22b,HZ19}
for related problems including PDEs. 

\item{(2)} \,
As mentioned at the end of Section \ref{sec:processes},
when the present non-Hermitian matrix-valued BM
starts from the null matrix, $M(0)=O$,
the statistics of eigenvalues and eigenvector-overlaps
at each time $t >0$ are identified with those in
the Ginibre ensemble with variance $t/N$.
It implies that this time-evolution of the system
can be regarded as a simple dilatation 
centered at the origin in $\C$ by the spacial factor $\sqrt{t}$
of the Ginibre statistics, which we write as
$(\Xi^{\rm Ginibre}(t, \cdot), \cO^{\rm Ginibre}(t, \cdot))_{t \geq 0}$.
In the present paper, we have studied the process
$(M(t))_{t \geq 0}$ starting from an arbitrary matrix $M(0)$.
If the supports of $\Xi(0, \cdot)$ and $\Theta(0, \cdot)$ 
are bounded, then we expect the convergence
of the processes $(\Xi(t, \cdot), \cO(t, \cdot))_{t \in [T, \infty)}$
to the simple `Ginibre process' 
$(\Xi^{\rm Ginibre}(t, \cdot), \cO^{\rm Ginibre}(t, \cdot))_{t \in [T, \infty)}$, 
as $T \to \infty$. 
Mathematical justification of such convergence will be
a future problem. 
In the $N \to \infty$ limit, the Ginibre ensemble exhibits
the celebrated \textit{circular law} 
which is universal in a wide class of non-Hermitian random-matrix
ensembles \cite{Bai97,Gir84,TV08}. 
The above convergence in $T \to \infty$
suggests a \textit{dynamical universality}
regardless of any details of bounded initial matrices.

\item{(3)} \,
Ginibre studied the statistical ensembles of
Gaussian and non-Hermitian random matrices 
with not only complex-valued entries,
but also with real-valued and quaternion-valued entries \cite{Gin65}.
The present non-Hermitian matrix-valued BM, $(M(t))_{t \geq 0}$, 
provides a dynamical extension of the 
\textit{complex} Ginibre ensemble.
A restriction to real-valued  and an extension to quaternion-valued BMs for the
entries of non-Hermitian matrix-valued BM 
should be studied in future. 
In the Hermitian random-matrix ensembles,
the parameter $\beta$ has been introduced, which takes
special values; $\beta=1$ for the real-valued, 
$\beta=2$ for the complex-valued,
and $\beta=4$ for the quaternion-valued entries.
Dyson derived the following system of SDEs for
the eigenvalue process
$(\bLambda^{\rm H}(t))_{t \geq 0}$ 
in the cases with $\beta=1,2$ and 4
for the Hermitian matrix-valued BM, $(M^{\rm H}(t))_{t \geq 0}$
\cite{Dys62,Kat_Springer},
\begin{equation}
d \Lambda_j^{\rm H}(t)= d B_j(t) + \frac{\beta}{2}
\sum_{1 \leq k \leq N, k \not=j}
\frac{dt}{\Lambda_j^{\rm H}(t)-\Lambda_k^{\rm H}(t)},
\quad 1 \leq j \leq N, \, t \geq 0.
\label{eq:Dyson}
\end{equation}
The interacting particle systems following 
\eqref{eq:Dyson} with general $\beta>0$
is now called the \textit{$\beta$-Dyson BM}
and extensively studied (see, for instance, \cite{For10}). 
Further considerations on the 
$\beta$-Dyson BM are found in \cite{HIM23,KK21}.
Is it possible to introduce such a relevant parameter $\beta$ into 
the present system of SDEs for
the eigenvalues \eqref{eq:SDE_Lambda1} 
and eigenvector-overlaps \eqref{eq:SDEofdOjk}, 
and into the SPDEs 
for the regularized FK-determinant 
random-fields \eqref{eq:SDE_FK}--\eqref{eq:SDE_psi1}?
See \cite{MT23} for the \textit{static} non-Hermitian 
$\beta$-ensemble. 

\item{(4)} \,
As mentioned in Section \ref{sec:processes}, 
if we consider the Hermitian matrix-valued process
$(M^{\rm H}(t))_{t \geq 0}$ 
the eigenvector-overlap matrix is identically equal to the
identity matrix; $\cO^{\rm H}(t) \equiv I$, $t \geq 0$. 
That is,
\begin{equation}
d \cO^{\rm H}_{jk}(t)=0, 
\quad 1 \leq j, k \leq N, \, t \geq 0.
\label{eq:Dyson2}
\end{equation}
In order to understand such significant difference
between the system of SDEs in the non-Hermitian case,
$(M(t))_{t \geq 0}$, studied in this paper
and the system of \eqref{eq:Dyson} and \eqref{eq:Dyson2}
in the Hermitian case, $(M^{\rm H}(t))_{t \geq 0}$, 
intermediate processes between
these two cases should be studied.
Recently, one of the present authors \cite{Yab20}
introduced the matrix-valued 
stochastic process 
with a parameter
$\tau \in [-1,1]$, 
$(M^{(\tau)}(t))_{t \geq 0}$, which interpolates
the two matrix-valued BMs. 
This process can be regarded as
a dynamical extension of the Girko ensembles
of random matrices \cite{Gir85,KS11,SCSS88}.
There the cross-variations of the increments
of the elements of $M^{(\tau)}(t)$ are given by
\begin{align*}
\bra d M^{(\tau)}_{jk}, d M^{(\tau)}_{\ell m} \ket_t
&= \tau \delta_{j m} \delta_{k \ell} \frac{dt}{N},
\nonumber\\
\bra d M^{(\tau)}_{jk}, d \overline{M^{(\tau)}_{\ell m}} \ket_t
&= \delta_{j \ell} \delta_{k m} \frac{dt}{N}, 
\quad 1 \leq j, k, \ell, m \leq N, \, t \geq 0.
\end{align*}
Hence, we can see that
$(M^{(1)}(t))_{t \geq 0} \law= (M^{\rm H}(t))_{t \geq 0}$
induces the eigenvalue process of the Dyson model 
\eqref{eq:Dyson} with $\beta=2$
on the real axis,
and that
$(M^{(0)}(t))_{t \geq 0} \law= (M(t))_{t \geq 0}$. 
We see that $(M(t)^{(-1)})_{t \geq 0}$ 
has a pure-imaginary-valued eigenvalue process 
which exhibits the Dyson model with $\beta=2$
on the imaginary axis \cite{Yab20}. 
SDEs for the eigenvalue process and 
eigenvector-overlap process,
and SPDEs for the regularized FK-determinant random fields will be studied 
for $(M^{(\tau)}(t))_{t \geq 0}$ with $\tau \in [-1, 1]$.

\end{description}
\vskip 0.3cm
\noindent{\bf Acknowledgements} \,
The present authors would like to thank 
Jacek Ma{\l}ecki, 
Tomoyuki Shirai, 
Govind Menon,
Horng-Tzer Yau,
Colin McSwiggen, 
Yan V. Fyodorov, 
Ping Zhong,
Hiroshi Kawabi, 
Shinji Koshida, 
Ryosuke Sato, 
Noriyoshi Sakuma, 
Makoto Nakashima, 
Kohei Noda, 
Yuya Tanaka, 
Saori Morimoto, 
and Ayana Ezoe 
for useful discussion on 
the present work.
This work was supported by the Research Institute for Mathematical Sciences,
an International Joint Usage/Research Center located in Kyoto University.
It was also supported (in part)
by funding from Fukuoka University
(Grant No.215001).
SE was supported by
JSPS KAKENHI Grant Numbers 
JP17K14206,
JP22H04942,
and
JP23K03158.
MK was supported by
JSPS KAKENHI Grant Numbers
JP16H06338,
JP18H01124,
JP19K03674,
JP21H04432, 
JP22H05105,
JP23H01077, 
JP23K25774, 
and
JP24K06888.
SY was supported by JSPS KAKENHI Grant Number
JP24K16940.


\end{document}